\documentclass[11pt,a4j]{article}
\usepackage{amsmath,amsthm,amssymb,bm}
\usepackage[authoryear,round]{natbib}
\usepackage{graphicx}
\usepackage{algorithm}
\usepackage{algorithmic}
\usepackage{color, xcolor}
\usepackage{breqn}

\usepackage{indentfirst}
\usepackage{url}
\def\X{{\bm{X}}}
\def\x{{\bm{x}}}
\def\F{{\bm{F}}}
\def\eps{{\bm{\varepsilon}}}
\def\0{{\bm{0}}}
\def\q{{q}}
\def\f{{f}}
\def\bfmu{\bm \mu}
\def\factanal{{\tt factanal} }
\def\R{{\tt R}}

\definecolor{darkgreen}{rgb}{0.0, 0.5, 0.0}

\newtheoremstyle{italicTheorem}% name
  {3pt}% Space above
  {3pt}% Space below
  {\itshape}% Body font
  {}% Indent amount
  {\itshape}% Theorem head font
  {.}% Punctuation after theorem head
  {.5em}% Space after theorem head
  {}% Theorem head spec (can be left empty, meaning ‘normal’)

\newcommand{\argmin}{\mathop{\rm arg~min}\limits}
\theoremstyle{plain}
\newtheorem{thm}{Theorem}
\newtheorem{pro}{Proposition}
\newtheorem*{lem*}{Lemma}
\newtheorem*{thm*}{Theorem}
\newtheorem{ex}{Example}
\theoremstyle{italicTheorem}
\newtheorem*{rem*}{Remark}
\newtheorem{rem}{Remark}
\theoremstyle{definition}
\newtheorem{dfn}{Definition}
\newtheorem*{dfn*}{Definition}

\allowdisplaybreaks[1]

\begin{document}

{
\begin{center}
  %\textbf{\Large Computer algebraic approach to maximum likelihood factor analysis
  \textbf{\Large Algebraic approach to maximum likelihood factor analysis    
}
\end{center}
\begin{center}
\large {Ryoya Fukasaku$^{1}$, Kei Hirose$^{2}$, Yutaro Kabata$^{3}$ and Keisuke Teramoto$^{4}$}
\end{center}

\begin{flushleft}
{\footnotesize
$^1$ Faculty of Mathematics, Kyushu University, 744 Motooka, Nishi-ku, Fukuoka 819-0395, Japan \\

$^2$ Institute of Mathematics for Industry, Kyushu University, 744 Motooka, Nishi-ku, Fukuoka 819-0395, Japan \\

$^3$ School of Information and Data Sciences, Nagasaki University, Bunkyocho 1-14, Nagasaki 852-8131, Japan \\

$^4$ Department of Mathematics, Yamaguchi University, Yoshida 1677-1, Yamaguchi 753-8512, Japan \\

\vspace{1.2mm}

}
{\it {\small E-mail: hirose@imi.kyushu-u.ac.jp, fukasaku@math.kyushu-u.ac.jp, kabata@nagasaki-u.ac.jp, kteramoto@yamaguchi-u.ac.jp}}	
\end{flushleft}

\vspace{1.5mm}

\begin{abstract}
In exploratory factor analysis, model parameters are usually estimated by maximum likelihood method.  The maximum likelihood estimate is obtained by solving a complicated multivariate algebraic equation.  Since the solution to the equation is usually intractable, it is typically computed with continuous optimization methods, such as Newton-Raphson methods.  With this procedure, however, the solution is inevitably dependent on the estimation algorithm and initial value since the log-likelihood function is highly non-concave. Particularly, the estimates of unique variances can result in zero or negative, referred to as improper solutions; in this case, the maximum likelihood estimate can be severely unstable.  To delve into the issue of the instability of the maximum likelihood estimate, we compute exact solutions to the multivariate algebraic equation by using algebraic computations. We provide a computationally efficient algorithm based on the algebraic computations specifically optimized for maximum likelihood factor analysis.  To be specific, Gr\"oebner basis and cylindrical decomposition are employed, powerful tools for solving the multivariate algebraic equation.   Our proposed procedure produces all exact solutions to the algebraic equation; therefore, these solutions are independent of the initial value and estimation algorithm. We conduct Monte Carlo simulations to investigate the characteristics of the maximum likelihood solutions.

 \end{abstract}

 \section{Introduction}
Factor analysis is a practical tool for investigating the covariance structure of observed variables by constructing a small number of latent variables referred to as common factors. Factor analysis has been used in the social and behavioral sciences since its proposal more than 100 years ago \citep{Spearman.1904}, and it has been recently applied to a wide variety of fields, including marketing, life sciences, materials sciences, and energy sciences \citep{Lin.2019,Shkeer.2019,Shurrab.2019,Kartechina.2020,Vilkaite-Vaitone.2022}.

In exploratory factor analysis, the model parameter is often estimated by the maximum likelihood method under the assumption that the observed variables follow the multivariate-normal distribution.  The maximum likelihood estimate is obtained by solving a complicated multivariate algebraic equation.   The solution to the equation is usually intractable, and thus it is typically computed using continuous optimization techniques.  One of the most common approaches is to use the gradient of the log-likelihood function, such as the Newton-Raphson method or Quasi-Newton method.  For given unique variances, the solution to the loading matrix turns out to be the eigenvalue/eigenvector problem for the standardized sample covariance matrix, which is similar to the problem of the principal component analysis \citep{Tipping.1999}.  Substituting the loading matrix into the log-likelihood function, it is regarded as a function of unique variances, and these unique variances are optimized by the Newton-Raphson method or Quasi-Newton method \citep{Joreskog.1967,Jennrich.1969,lawley1971}.  The above algorithm has been widely-used for several decades; indeed, \factanal function, one of the most popular functions that implements the maximum likelihood factor analysis in \R, is based on \citeauthor{lawley1971}'s \citeyearpar{lawley1971} algorithm.  Another algorithm is the expectaion-maximization (EM) algorithm \citep{RUBIN.1982}, in which the unique factors are regarded as latent variables.  The complete-data log-likelihood function corresponds to the log-likelihood function in multivariate regression; thus, the EM algorithm can be applied to the penalized likelihood approach, such as the lasso and minimax concave penalty \citep{Choi.2010,Hirose.2015}, allowing to achieve high-dimensional sparse estimation.   

Although the above-mentioned algorithms often converge, the convergence value is not always optimal due to the non-concavity of the log-likelihood function. Figure \ref{fig:discrepancy3} shows the convergence values of the discrepancy function that corresponds to negative log-likelihood function and maximum likelihood estimates of unique variance for second observed variable, say $\hat{\psi}_2$. These convergence values are obtained by applying estimation algorithms to an artificial dataset, generated from simulation model {\rm S2} that will be described in Section 5. We employ three estimation algorithms: (i) the \factanal function in \R \ (i.e., \citeauthor{lawley1971}'s \citeyearpar{lawley1971} algorithm), (ii) the Quasi-Newton method based on \citet{Jennrich.1969}, and (iii) the EM algorithm \citep{RUBIN.1982}. 100 different initial values for unique variances are generated from a uniform distribution $U(0,1)$ to investigate whether the convergence value depends on the initial value.  It should be noted that we use the same sets of 100 initial values for these three algorithms.

\begin{figure}[!t]
\centering
\includegraphics[width=\textwidth, bb=0 0 804 343]{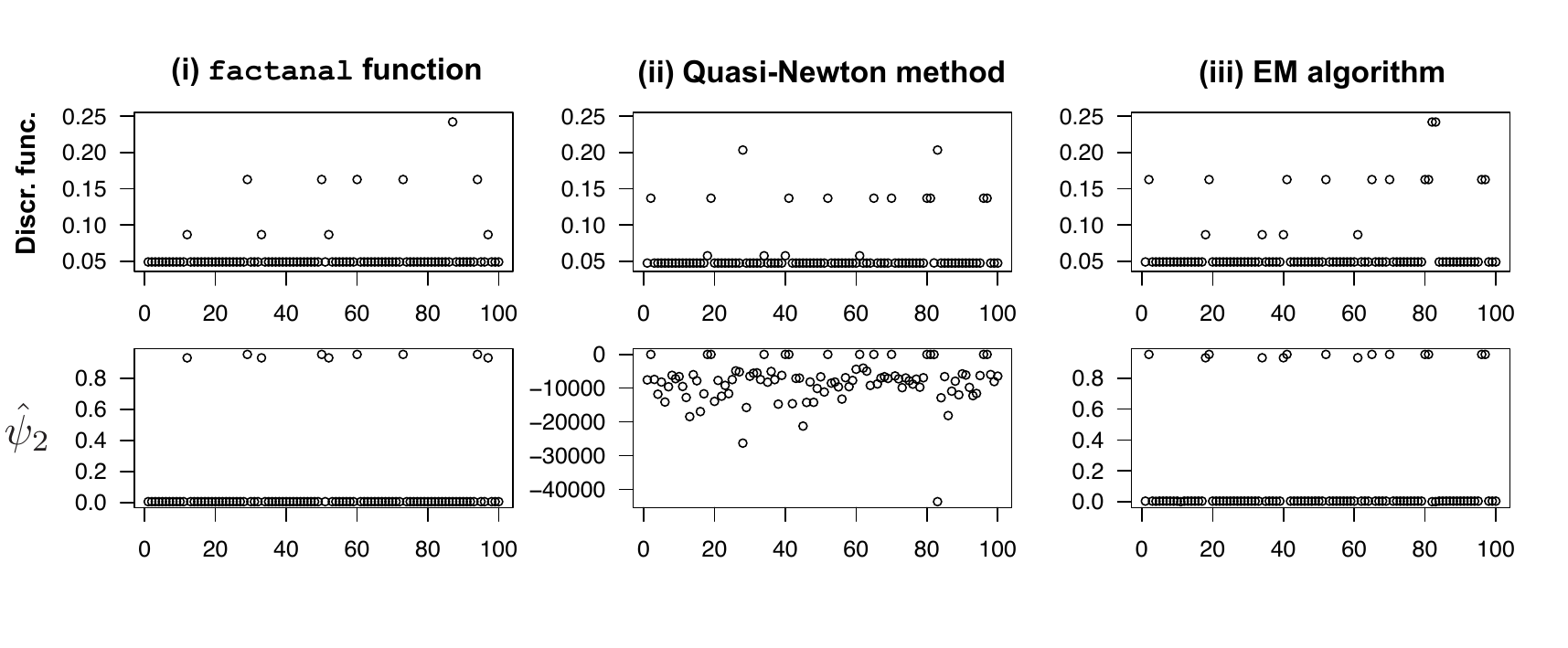}
\caption{Convergence values of the discrepancy function and maximum likelihood estimates of unique variance for second observed variable, $\hat{\psi}_2$.  These convergence values are obtained by applying three estimation algorithms to an artificial dataset with 100 different initial values of unique variances, and the horizontal axis indicates the indices for these 100 initial values.}
\label{fig:discrepancy3}
\end{figure}

The results indicate that the convergence values are dependent on the initial values and estimation algorithms. In particular, the estimate of unique variance, $\hat{\psi}_2$, is sensitive to the initial values with the Quasi-Newton method (i.e., \citeauthor{Jennrich.1969}'s \citeyearpar{Jennrich.1969} algorithm).  We observe that improper solutions are obtained; that is, the estimates of unique variances turn out to be zero or negative. The Quasi-Newton method allows the estimates of unique variances to take zero or negative values; thus, the $\hat{\psi}_2$ becomes exceptionally unstable.  Meanwhile, the \factanal function typically outputs parameters such that all of the unique variances are greater than some small threshold values, such as 0.005.  The EM algorithm also does not produce negative estimates of unique variances due to the characteristics of the algorithm \citep{Adachi.2013}. Although the \factanal function and EM algorithm do not yield negative unique variances, their results are dependent on the initial values.    Consequently, it is quite difficult to evaluate whether an optimal value of the parameter is obtained with these estimation algorithms.

As shown in the above example, the instability of the numerical solutions is often observed when we obtain an improper solution in our experience.  To our knowledge, a theoretical approach to elucidate the improper solutions has not yet been provided due to the complexity of the log-likelihood function; thus, the improper solutions remain to be a challenging problem. Instead of theoretical approaches, numerical approaches have been employed to investigate the characteristics of improper solutions for the last five decades (e.g., \citealp{Joreskog.1967,Driel.1978,Sato.1987,Kano.1998,Krijnen.1998,Hayashi.2014}).  However, with the numerical approaches, the result is highly dependent on the estimation algorithm and its initial values, as shown in Figure \ref{fig:discrepancy3}. 
%The fundamental issue is that the numerical approach is based on a continuous optimization problem
%obtained by solving a continuous optimization problem; through this procedure, the solution is inevitably dependent on the initial value and estimation algorithm.  

To address this issue, it is essential to compute {\it exact} solutions to the multivariate algebraic equation.  In this study, we employ an algebraic approach based on computational algebra to obtain the exact solutions.  Computational algebra can solve complicated multivariate algebraic equations by using the theory of Gr\"oebner bases. Several software that implement computational algebra, such as Magma, Maple, and Mathematica, have been developed considerably in the last two decades.  Thanks to the development of the algorithm and the rapid progress in computational technology, several real problems have been resolved by computational algebra (e.g., \citealp{Laubenbacher.2009,Hinkelmann.2011,Veliz-Cuba.2014}).   

In maximum likelihood factor analysis, we observe that the original multivariate algebraic equation requires a lot of computer resources from an algebraic viewpoint, such as inverse covariance matrix computation.  To reduce the heavy computational loads, we introduce a computationally efficient algorithm  specifically optimized for maximum likelihood factor analysis. In our algebraic algorithm, we employ the theory of Gr\"oebner bases to get simplified sub-problems for the algebraic equation. After getting the simplified sub-problems, we compute all exact solutions to the algebraic equation. The maximum likelihood estimate can be obtained by selecting a solution that maximizes the likelihood function.  The solution is independent of the initial value and estimation algorithm because our algebraic algorithm produces all exact solutions.

The algebraic approach to factor analysis has been explored by several earlier studies, as highlighted by  \citet{Drton2007,Drton.2010,Ardiyansyah.2022,Drton.2023}. In particular, \citet{Drton2007} employs computational algebra to study model invariants in factor analysis, and the final section spotlights algebraic study for maximum likelihood estimation and singularities as the next steps. These earlier works mostly focus on perspectives from mathematics. This study, however, aims to develop an efficient algorithm for computing the exact solution, enabling us to delve into the characteristics of the maximum likelihood solutions.  Specifically, an extensive Monte Carlo simulation is conducted to provide a  detailed analysis of these characteristics.  We observe that sometimes the dimension of the solution space can be more than zero.  Furthermore, the maximum likelihood estimate does not exist in some cases even if the parameter space of unique variances includes negative values. Such a discussion is impossible with existing algorithms.  Consequently, the proposed algorithm provides not only theoretical insights but also numerical analysis that significantly advance elucidating properties of maximum likelihood solutions, including the characterization of improper solutions, thereby contributing to a deeper understanding from a practical viewpoint. 

The remainder of this paper is organized as follows: Section 2 briefly reviews the maximum likelihood factor analysis. An algebraic equation to solve the maximum likelihood solution is also discussed. In Section 3, we introduce a novel algorithm for computing the exact maximum likelihood solutions via computational algebra.  Section 4 presents numerical results for artificial datasets. Some concluding remarks are given in Section 5.

%%%%%%%%%%%%%%%%%%%%%%%%%%%%%%%%%%%%%%%%%%%%%%%%%%%%%%%%%%%%%%%%%%%%%%%%%%%%%%%%%%%%%%%%%%%%%%%%%%%%%%%%%%%%%%%%%%%%%%%%%%%%%%%%%%%%%%%%%%%%%%%%%%%%%%%%%%%%%%%
\section{Maximum likelihood factor analysis}
\subsection{Maximum likelihood estimation}
Let $\X=(X_1,\cdots,X_p)^\top$ be a $p$-dimensional observed random vector. The factor analysis model is
\begin{equation*}
\X =\bfmu + L\F + \eps  , \label{model1}
\end{equation*}
where $\bfmu$ is a mean vector, $L =(\ell_{ij})$ is a $p \times k$  matrix of factor loadings, and $\F = (F_1,\cdots,F_k)^\top$ and $\eps  = (\varepsilon_1,\cdots, \varepsilon_p)^\top$ are unobservable random vectors. Here, $A^\top$ denotes the transpose of a matrix $A$. The elements of \mbox{\boldmath{$F$}}  and \mbox{\boldmath{$\varepsilon$}}  are referred to as common factors and unique factors, respectively. It is assumed that $\F \sim N_k(\0, I_k)$ and $\eps \sim N_p(\0, \Psi)$ and $\F$ and $\eps$ are independent. Here,  $I_k$ is an identity matrix of order $k$ and $\Psi $ is a $p \times p$ diagonal matrix with $i$-th diagonal element being $\psi_i$, referred to as unique variance.  Under these assumptions, the observed vector $\X$ follows multivariate-normal distribution;  $\X \sim N_p(\bfmu, \Sigma)$ with $\Sigma=LL^\top+\Psi$. When the observed variables are standardized, $i$-th diagonal element of $ LL^\top$ is called communality, which measures the percent of variance in $x_i$ explained by common factors. 

It is well-known that factor loadings have a rotational indeterminacy since both $L $ and  $LT$ generate the same covariance matrix $\Sigma$, where $T$ is an arbitrary orthogonal matrix.  To obtain an identifiable model, we often use a restriction that $L^\top\Psi^{-1}L$ is a diagonal matrix or the upper triangular elements of $L$ are fixed by 0.

Suppose that we have a random sample of $N$ observations $\x_1,\cdots,\x_N$ from the $p$-dimensional normal population $N_p(\bfmu ,\Sigma)$. The maximum likelihood estimate of factor loadings and unique variances is obtained by maximizing the log-likelihood function
\begin{equation}
\log \f(\bm{x}_1,\dots,\bm{x}_N | L,\Psi  ) = -\frac{N}{2} \bigg\{ p\log(2\pi)+\log |\Sigma| + \mathrm{tr}(\Sigma^{-1} S ) \bigg\}, \label{taisuuyuudo}	
\end{equation}
or equivalently, minimizing the discrepancy function:
\begin{equation}
\q( L,\Psi ) = \log|\Sigma|+ \mathrm{tr}(\Sigma^{-1}S) -\log |S|-p, \label{q}
\end{equation}
where $S=(s_{ij})$ is the sample variance-covariance matrix,
\begin{equation*}
S=\frac{1}{N} \sum_{n=1}^N(\x_n-\bar{\x})(\x_n-\bar{\x})^\top, \quad \bar{\bm{x}} = \frac{1}{n}\sum_{i=1}^n\bm{x}_i.
\end{equation*}
We assume that $S$ is non-singular to ensure the existence of $S^{-1}$, which plays an important role in constructing an algorithm to find a maximum likelihood solution.  This assumption can be violated when the number of variables exceeds the sample size. In such a case, we may employ the ridge regularization on $S$; that is, we use $S_{\lambda}:=S+\lambda I_p$ with $\lambda>0$ \citep{Chan.2008}.

\subsection{Numerical algorithms}
We briefly review the algorithms for computing the maximum likelihood estimates \citep{Joreskog.1967,lawley1971}. The candidates of the maximum likelihood estimates, say $(\hat{L},\hat{\Psi})$, are given as the solutions of 
\begin{equation}
\frac{\partial \q(L ,\Psi ) }{ \partial L}  = 0,~~ \frac{\partial \q(L,\Psi) }{ \partial \Psi}  = 0. \label{eq:diff}	
\end{equation}
The solution to the above equation is given by 
\begin{align}
  \begin{split}
\Sigma^{-1}(\Sigma-S)\Sigma^{-1}L &=0,\\
\rm{diag}(\Sigma^{-1}(\Sigma-S)\Sigma^{-1}) &=0.
  \label{eq:diff0}
  \end{split}
\end{align}

Since the solutions are not expressed in a closed form, some iterative algorithm is required. Typically,  we use an algorithm based on the gradient of the log-likelihood function, such as the Newton-Raphson method or Quasi-Newton method.  When $\Psi$ is not singular, Eq. \eqref{eq:diff} results in the following expressions  \citep{lawley1971}:
\begin{align}
(\Psi^{-\frac{1}{2}} S \Psi^{-\frac{1}{2}})  \Psi^{-\frac{1}{2}} L&= \Psi^{-\frac{1}{2}} L(I_k +L^\top\Psi^{-1}L), \label{eq:diff1-1} \\
\Psi&=\mathrm{diag}(S-LL^\top). \label{eq:diff1-2} 	
\end{align}
Let $L^\top\Psi^{-1}L=\Delta $.  The rotational indeterminacy is resolved when $\Delta$ is a diagonal matrix (e.g., \citealp{lawley1971}); in this case, \eqref{eq:diff1-1} is an eigenvalue/eigenvector problem.  Let $(\theta_1,\bm{e}_1),\cdots,(\theta_p,\bm{e}_p)$ be the eigenvalues and eigenvectors of $\Psi^{-\frac{1}{2}} S \Psi^{-\frac{1}{2}}$, respectively, where the eigenvalues are rearranged in a decreasing order; that is, $\theta_1 \geq \theta_2 \geq \cdots \geq \theta_p$.  Let $P =(\bm{e}_1,\cdots ,\bm{e}_k)$ and $\Theta =\mathrm{diag}(\theta_1,\theta_2, \cdots , \theta_k)$.  Eq. \eqref{eq:diff1-1} implies
\begin{equation}
P (\Theta -I_k )^{\frac{1}{2}}=\Psi^{-\frac{1}{2}}L \label{PTheta-I2-1}
\end{equation}
or equivalently,
\begin{align}
L&= \Psi^{\frac{1}{2}} P (\Theta -I_k )^{\frac{1}{2}}.  \label{Lsuiteirou2-1}	
\end{align}
Thus, for given $\Psi$, $L$ can be calculated using \eqref{Lsuiteirou2-1}.  For given $L$, one can update $\Psi$ with \eqref{eq:diff1-2}.  However, the update with  \eqref{eq:diff1-2} can be slow when some of the diagonal elements of $\hat{\Psi}$ are small.  Thus, the unique variances are often updated by using the Newton-Raphson method or quasi-Newton method \citep{Joreskog.1967,lawley1971}.  The \factanal function in \R \ uses the {\tt optim} function to update the unique variances.  

To summarize, the estimation algorithm in \citet{Joreskog.1967} is given as follows:
\begin{enumerate}
\item Set an initial value of $\Psi$.
\item For given $\Psi$, Eq. \eqref{Lsuiteirou2-1} implies that $L$ is a function of $\Psi$; that is, $L$ can be written as $L=L(\Psi)$. Substituting $L(\Psi)$ into the discrepancy function \eqref{q}, it is a function of $\Psi$.  The elements of $\Psi$ are then minimized by Newton-Raphson method or quasi-Newton method. 
\end{enumerate}

\subsection{An algebraic equation of the maximum likelihood solution}
%Although the maximum likelihood solution is usually obtained by using an algorithm in \citet{lawley1971}, the estimate of unique variances can sometimes be exactly zero or negative, referred to as improper solutions.  In this case, the solution can be numerically unstable as shown in Figure \ref{fig:discrepancy3}.  Thus, 
In this study, we will employ a computational algebraic approach to get an exact solution to algebraic equations \eqref{eq:diff0}.  However, Eq. \eqref{eq:diff0} may not be directly used as it includes an inverse matrix $\Sigma^{-1}$ that consists of significantly complicated rational functions.  Consequently, simpler algebraic equations are required.  One may use \eqref{eq:diff1-1} -- \eqref{eq:diff1-2} as they do not include $\Sigma^{-1}$, but Eq. \eqref{eq:diff1-1} includes $\Psi^{-1}$, which is unavailable when $\Psi$ is singular.  As will be demonstrated in the numerical example in Section \ref{sec:MCS}, there exist many solutions whose unique variances are exactly zeros.  Thus, it is essential to derive an algebraic equation that does not include complicated rational functions and also satisfies \eqref{eq:diff0}, even when $\Psi$ is singular.

In the following theorem, we provide a necessary condition for \eqref{eq:diff0} that does not depend on  rational functions and holds even when $\Psi$ is singular.  Therefore, all solutions of \eqref{eq:diff0} must satisfy that necessary condition. 
\begin{thm}\label{we-want-solve}
Whether $\Psi$ is singular or not, if \eqref{eq:diff0} is satisfied, then we have
  \begin{align}
    L&=(LL^\top + \Psi) S^{-1}L, \label{formula-L} \\
    \Psi&=\mathrm{diag}(S-LL^\top). \tag{\ref{eq:diff1-2}}
  \end{align}
  \begin{proof}
It is shown that Eq. \eqref{formula-L} is equivalent to Eq. \eqref{eq:diff0}, since we have
    \begin{align*}
      \Sigma^{-1}(\Sigma-S)\Sigma^{-1}L =0
      &\Longleftrightarrow
      (\Sigma S^{-1} \Sigma) \Sigma^{-1}(\Sigma-S)\Sigma^{-1}L =0
      \\
      &\Longleftrightarrow
      \Sigma S^{-1}L=L
      \\
      &\Longleftrightarrow
      (LL^\top + \Psi) S^{-1}L=L.
    \end{align*}
Therefore, what we need to prove is that Eq. \eqref{eq:diff0} implies Eq. \eqref{eq:diff1-2}.  First, we have
    \begin{align*}
      \Sigma - S = (LL^\top + \Psi) \Sigma^{-1} (\Sigma - S) \Sigma^{-1} (LL^\top + \Psi) 
    \end{align*}
    by $\Sigma = LL^\top + \Psi$. Hence, the upper part of \eqref{eq:diff0} implies that
    \begin{align*}
      \Sigma - S = \Psi \Sigma^{-1} (\Sigma - S) \Sigma^{-1} \Psi.
    \end{align*}
To handle the case where $\Psi$ is singular, suppose that $\psi_1, \ldots, \psi_q$ are not equal to zero and $\psi_{q+1}, \ldots, \psi_p$ are equal to zero, where $q \leq p$. Let $J = \{ 1, \ldots, q\}, K = \{q+1,\ldots,p\}$, and let
    \begin{align*}
      T = \mathrm{diag}\left(\frac{1}{\psi_1}, \ldots, \frac{1}{\psi_q}, 1, \ldots, 1\right), ~~ I_{q \leq p} = \mathrm{diag}\left(\overbrace{1, \ldots, 1}^{q}, \overbrace{0, \ldots, 0}^{p-q}\right).
    \end{align*}
Let $M_{KJ}$ denote a submatrix of a $M$ from which the rows corresponding to the indices in $K$ and the columns corresponding to the indices in $J$ are removed, respectively.  Since $T$ is not singular,    
    \begin{align*}
      \Sigma - S = \Psi \Sigma^{-1} (\Sigma - S) \Sigma^{-1} \Psi
      &\Longleftrightarrow
      T(\Sigma - S) T = T \Psi \Sigma^{-1} (\Sigma - S) \Sigma^{-1} \Psi T
      \\
      &\Longleftrightarrow
      T(\Sigma - S) T = I_{q \leq p} \Sigma^{-1} (\Sigma - S) \Sigma^{-1} I_{q \leq p}
      \\
      &\Longleftrightarrow
      T(\Sigma - S) T =
      \begin{pmatrix}
        (\Sigma^{-1} (\Sigma - S) \Sigma^{-1})_{KK} & 0
        \\
        0 & 0
      \end{pmatrix}.
    \end{align*}
    As $(T (\Sigma - S) T)_{KK} = (\Sigma^{-1} (\Sigma - S) \Sigma^{-1})_{KK}$ holds, the lower part of \eqref{eq:diff0} implies that
    \begin{align*}
      \mathrm{diag}((T (\Sigma - S) T)_{KK}) = 0.
    \end{align*}
Hence, $\Psi_{KK} = \mathrm{diag}((S - LL^\top)_{KK})$ holds because $T$ is diagonal. Since
    \begin{align*}
      (\Sigma - S)_{JJ} = I_{p-q} (\Sigma - S)_{JJ} I_{p-q} = (T (\Sigma - S) T)_{JJ} = 0
    \end{align*}
    holds, we also have $s_{ii} - \sum_{j=1}^k \ell_{ij}^2 = 0 = \psi_i$ for each $i \in K$. Consequently, we obtain Eq. \eqref{eq:diff1-2} without an assumption that $\Psi$ is not singular.
  \end{proof}
\end{thm}

%\begin{rem}
%Eq. \eqref{formula-L} is much more useful than Eq. \eqref{eq:diff1-1} in terms of computational efficiency in the algebraic approach. In fact, comparing Eq. \eqref{formula-L} and the following equation
%\begin{align}
%    L &= (S - LL^\top) \Psi^{-1} L, \label{eq:diff1-2-2}
%\end{align}
%which is equivalent to \eqref{eq:diff1-1},  \eqref{eq:diff1-2-2} contains rational functions while Eq. \eqref{formula-L} does not. In particular, since each $\psi_i$ is represented as a polynomial in variables $\ell_{ij}$ from Eq. \eqref{eq:diff1-2}, we obtain an significantly complicated rational functions by substituting \eqref{eq:diff1-2} into \eqref{eq:diff1-2-2}.
%\end{rem}

\begin{rem}
As an alternative algorithm to obtain the maximum likelihood estimates given by \citet{lawley1971}, \citet{Jennrich.1969} used the eigenvalue/eigenvector problem of $S^{-\frac{1}{2}} \Psi S^{-\frac{1}{2}}$ instead of $\Psi^{-\frac{1}{2}} S \Psi^{-\frac{1}{2}}$ by using the following equation: 
\begin{align}
(S^{-\frac{1}{2}} \Psi S^{-\frac{1}{2}})  S^{-\frac{1}{2}} L&= S^{-\frac{1}{2}} L(I_k -L^\top S^{-1}L),\label{Ldiff2}
\end{align}
which is equivalent to Eq. \eqref{formula-L}.  Hence, \citeauthor{Jennrich.1969}'s \citeyearpar{Jennrich.1969} algorithm can produce zero or negative unique variance estimates. 
\end{rem}

%%%%%%%%%%%%%%%%%%%%%%%%%%%%%%%%%%%%%%%%%%%%%%%%%%%%%%%%%%%%%%%%%%%%%%%%%%%%%%%%%%%%%%%%%%%%%%%%%%%%%%%%%%%%%%%%%%%%%%%%%%%%%%%%%%%%%%%%%%%%%%%%%%%%%%%%%%%%%%%
\section{Maximum likelihood solution via Gr\"oebner basis}

Theorem \ref{we-want-solve} implies all solutions to \eqref{eq:diff0} must satisfy \eqref{eq:diff1-2} and \eqref{formula-L}; therefore, we first solve all exact solutions to \eqref{eq:diff1-2} and \eqref{formula-L} to find the candidates of \eqref{eq:diff0}, and then obtain an optimal solution that maximizes the likelihood function.  In this section, we propose an algebraic approach to finding exact solutions to \eqref{eq:diff1-2} and \eqref{formula-L}.  Section \ref{mlsgb-1} reviews a theory of Gr\"obner bases for algebraic equations (see \citealp{CLO} for details). A simple illustrative example of the Gr\"oebner basis in maximum likelihood factor analysis is also presented.  We provide an essential idea for constructing an algorithm to find the exact maximum likelihood solution through the illustrative example. In Section \ref{mlsgb-2}, we introduce an algorithm for finding all candidates for the maximum likelihood solutions by using the theory of the Gr\"obner basis.  Particularly, we decompose an original algebraic equation into several sub-problems specifically to the maximum likelihood factor analysis, allowing us to get the maximum likelihood solution within a feasible computation time.  Furthermore, Section \ref{mlsgb-3} proposes an algorithm to obtain an optimal solution from the candidates for maximum likelihood solutions. Our proposed algorithm identifies a pattern of maximum likelihood solution among ``Proper solution", ``Improper solution", and ``No solutions." 

%%%%%%%%%%%%%%%%%%%%%%%%%%%%%%%%%%%%%%%%%%%%%%%%%%%%%%%%%%%%%%%%%%%%%%%%%%%%%%%%%%%%%%%%%%%%%%%%%%%%%%%%%%%%%%%%%%%%%%%%%%%%%%%%%%%%%%%%%%%%%%%%%%%%%%%%%%%%%%%
\subsection{Algebraic equations and Gr\"obner bases}\label{mlsgb-1}
\subsubsection{Affine varieties and algebraic equations}
Let $\mathbb{R}[\bm{z}]$ denote the set of all polynomials in $m$ variables $\bm{z} = (z_1, \ldots, z_m)$ with coefficients in the real number field $\mathbb{R}$. Recall that $\mathbb{R}[\bm{z}]$ is a commutative ring, which is called the polynomial ring in variables $\bm{z}$ with coefficients in $\mathbb{R}$. For details of fields and commutative rings, please refer to \ref{Fields}.  Let $f_1, \ldots, f_r$ be polynomials in the polynomial ring $\mathbb{R}[\bm{z}]$. We consider an algebraic equation which has the following general form:
\begin{align}
  \left\{
  \begin{array}{c}
    f_1(\bm{z}) = 0,
    \\
    \vdots
    \\
    f_r(\bm{z}) = 0.
  \end{array}
  \right.
  \label{AEGB-1}
\end{align}
Define
\begin{align*}
  \mathbb{V}_{W}(f_1, \ldots, f_r) = \{ \bm{z} \in W^m : f_1(\bm{z}) = \cdots = f_r(\bm{z}) = 0\} \subset W^m,
\end{align*}
where $W = \mathbb{R}$ or $\mathbb{C}$. The $ \mathbb{V}_{W}(f_1, \ldots, f_r)$ is referred to as affine variety defined by $f_1, \ldots, f_r$ in $W^{m}$. Note that the solutions of \eqref{AEGB-1} in $W^m$ form the affine variety $\mathbb{V}_{W}(f_1,\ldots, f_r)$.

Here, let us consider an illustrative example related to the algebraic equations \eqref{eq:diff1-2} and \eqref{formula-L}.

\begin{ex}\label{ex1}
  Let
  \begin{align*}
    S =
    \begin{pmatrix}
      1 & \frac{1}{2} & \frac{1}{3}
      \\
      \frac{1}{2} & 1 & \frac{2}{3}
      \\
      \frac{1}{3} & \frac{2}{3} & 1
    \end{pmatrix},
    ~~
    L =
    \begin{pmatrix}
      \ell_{11}
      \\
      \ell_{21}
      \\
      \ell_{31}
    \end{pmatrix},
    ~~
    \Psi =
    \begin{pmatrix}
      \psi_1 & 0 & 0
      \\
      0 & \psi_2 & 0
      \\
      0 & 0 & \psi_3
    \end{pmatrix}.
  \end{align*}
  We consider the algebraic equations \eqref{eq:diff1-2} and \eqref{formula-L} for the above $S, L, \Psi$, 
  that is,
  \begin{align}
  \left\{\begin{array}{l}
    0 = (\Psi - \mathrm{diag}(S - LL^\top))_{11},
    \\
    0 = (\Psi - \mathrm{diag}(S - LL^\top))_{22},
    \\
    0 = (\Psi - \mathrm{diag}(S - LL^\top))_{33},
    \\
    0 = (L - (LL^\top + \Psi)S^{-1}L)_{11},
    \\
    0 = (L - (LL^\top + \Psi)S^{-1}L)_{21},
    \\
    0 = (L - (LL^\top + \Psi)S^{-1}L)_{31},
  \end{array}\right.
  \label{exaegb-1}
  \end{align}
  We give polynomials $f_1, \ldots, f_6 \in \mathbb{R}[\bm{\psi}, \bm{\ell}] = \mathbb{R}[\psi_1, \psi_2, \psi_3, \ell_{11},\ell_{21},\ell_{31}]$ by
  \begin{align*}
      f_1 &= (\Psi - \mathrm{diag}(S - LL^\top))_{11} = \psi_1 - (1 - \ell_{11}^2),
      \\
      f_2 &= (\Psi - \mathrm{diag}(S - LL^\top))_{22} = \psi_2 - (1 - \ell_{21}^2),
      \\
      f_3 &= (\Psi - \mathrm{diag}(S - LL^\top))_{33} = \psi_3 - (1 - \ell_{31}^2),
      \\
      f_4 &= (L - (LL^\top + \Psi)S^{-1}L)_{11} = \ell_{11} - 
      \left(\begin{array}{l}
        \ell_{11} \left(- \frac{2}{3} \ell_{11} \ell_{21} + \frac{4}{3} (\ell_{11}^2 + \psi_1)\right)
        \\
        + \; \ell_{21} \left(\frac{32}{15} \ell_{11} \ell_{21} - \frac{6}{5} \ell_{11} \ell_{31} - \frac{2}{3} (\ell_{11}^2 + \psi_1)\right)
        \\
        + \; \ell_{31} \left(-\frac{6}{5} \ell_{11} \ell_{21} + \frac{9}{5}\ell_{11} \ell_{31}\right)
      \end{array}\right),
      \\
      f_5 &= (L - (LL^\top + \Psi)S^{-1}L)_{21} = \ell_{21} -
      \left(\begin{array}{l}
        \ell_{11} \left(\frac{4}{3} \ell_{11} \ell_{21} - \frac{2}{3} (\ell_{21}^2 + \psi_2)\right)
        \\
        + \; \ell_{21} \left(-\frac{2}{3} \ell_{11} \ell_{21} - \frac{6}{5} \ell_{21} \ell_{31} + \frac{32}{15} (\ell_{21}^2 + \psi_2)\right)
        \\
        + \; \ell_{31} \left(\frac{9}{5} \ell_{21} \ell_{31} - \frac{6}{5} (\ell_{21}^2 + \psi_2)\right)
      \end{array}\right),
      \\
      f_6 &= (L - (LL^\top + \Psi)S^{-1}L)_{31} = \ell_{31} -
      \left(\begin{array}{l}
        \ell_{11} \left(\frac{4}{3} \ell_{11} \ell_{31} - \frac{2}{3} \ell_{21} \ell_{31}\right)
        \\
        + \; \ell_{21} \left(-\frac{2}{3} \ell_{11} \ell_{31} + \frac{32}{15} \ell_{21} \ell_{31} - \frac{6}{5} (\ell_{31}^2 + \psi_3)\right)
        \\
        + \; \ell_{31} \left(-\frac{6}{5} \ell_{21} \ell_{31} + \frac{9}{5} (\ell_{31}^2 + \psi_3)\right)
        \end{array}\right).
  \end{align*}
  All solutions of the above equation forms the affine variety defined by $f_1,\ldots,f_6$. Indeed, there are seven real solutions to \eqref{exaegb-1}, and these solutions, $(\bm{\psi}, \bm{\ell}) = (\psi_1, \psi_2, \psi_3, \ell_{11}, \ell_{21}, \ell_{31}) \in \mathbb{R}^6$, form the following affine variety: 
  \begin{align*}
    \mathbb{V}_{\mathbb{R}}(f_1,\ldots,f_6)
    =
    \left\{\begin{array}{ll}
    \left(0, \frac{3}{4}, \frac{8}{9}, \pm 1, \pm \frac{1}{2}, \pm \frac{1}{3}\right),
    &
    \left(\frac{3}{4}, 0, \frac{5}{9}, \pm \frac{1}{2}, \pm 1, \pm \frac{2}{3} \right),
    \\
    \left(\frac{8}{9}, \frac{5}{9}, 0, \pm \frac{1}{3}, \pm \frac{2}{3}, \pm 1\right),
    &
    \left(1, 1, 1, 0, 0, 0\right)
    \end{array}\right\}.
  \end{align*}
\end{ex}

\subsubsection{Sums and saturations}
In Example \ref{ex1}, we obtain only seven solutions to the algebraic equations \eqref{exaegb-1}. For larger model, however, the number of solutions of \eqref{eq:diff1-2} and \eqref{formula-L} is considerably large and the solution space (i.e., the affine variety) is too complicated.  In some cases, there are an infinite number of solutions. Therefore, it is essential to enhance an efficient computation for obtaining all solutions. A key idea is to a decomposition of \eqref{eq:diff1-2} and \eqref{formula-L} into several simple sub-problems. This decomposition is related to the theory of polynomial ideals and Gr\"obner basis.

Let $J$ be the ideal generated by $f_1, . . . ,f_r$, that is,
\begin{align*}
  \mathcal{J} = \langle f_1, \ldots, f_r \rangle = \left\{\sum_{i=1}^r q_i f_i : q_i \in \mathbb{R}[\bm{z}]\right\} \subset \mathbb{R}[\bm{z}].
\end{align*}
For detail of the ideal in the commutative ring, please refer to \ref{Fields}.  The affine variety defined by $\mathcal{J}$ in $W^m$ is expressed as 
$$\mathbb{V}_{W}(\mathcal{J}) = \{ \bm{z} \in W^m : f(\bm{z}) = 0, \forall f \in \mathcal{J}\}.$$  
Since $\mathbb{V}_{W}(f_1, \ldots, f_r) = \mathbb{V}_{W}(\mathcal{J})$ holds, the solution space of \eqref{AEGB-1} is $\mathbb{V}_{W}(\mathcal{J})$.  

Since ideals are algebraic objects, there exist natural algebraic operations.  Essentially, such algebraic operations give us some simple sub-problems of the algebraic equations \eqref{eq:diff1-2} and \eqref{formula-L}. Let $\mathcal{K}$ be an ideal in $\mathbb{R}[\bm{z}]$. We note that any ideal in $\mathbb{R}[\bm{z}]$ is finitely generated by the Hilbert Basis Theorem (\citealp{CLO}; theorem 4, section 2.5), that is, there exists $h_1, \ldots, h_s \in \mathcal{K}$ such that $\mathcal{K} = \langle h_1, \ldots, h_s \rangle$.
%$h_1, \ldots, h_s \in \mathbb{R}[\bm{z}]$

Now, we define two algebraic operations. As the first algebraic operation for ideals, we define the sum of $\mathcal{J}$ and $\mathcal{K}$, say $\mathcal{J} + \mathcal{K}$, by
\begin{align*}
  \mathcal{J} + \mathcal{K} = \{ f + h : f \in \mathcal{J}, h \in \mathcal{K} \}.
\end{align*}
Note that any sum of ideals is an ideal in $\mathbb{R}[\bm{z}]$ (\citealp{CLO}; proposition 2, section 4.3). Moreover, the affine variety $\mathbb{V}_W(\mathcal{J} + \mathcal{K})$ coincides with the affine variety $\mathbb{V}_W(\mathcal{J}) \cap \mathbb{V}_W(\mathcal{K})$ (\citealp{CLO}; corollary 3, section 4.3), that is, $\mathbb{V}_W(\mathcal{J} + \mathcal{K})$ is none other than the solution space of the equation
\begin{align*}
    \eqref{AEGB-1}
    \mbox{ and }
    \left\{
    \begin{array}{c}
      h_1(\bm{z}) = 0,
      \\
      \vdots
      \\
      h_s(\bm{z}) = 0.
    \end{array}
    \right.  
\end{align*}

As the second algebraic operations, we define the saturation of $\mathcal{J}$ with respect to $\mathcal{K}$, say $\mathcal{J}:\mathcal{K}^{\infty}$, by
\begin{align*}
  \mathcal{J}:\mathcal{K}^{\infty} = \{ f \in \mathbb{R}[\bm{z}] : \mbox{for any } h \in \mathcal{K} \mbox{ there exists } t \in \mathbb{Z}_{\geq 0}  \mbox{ such that } fh^t \in \mathcal{J} \},
\end{align*}
where $\mathbb{Z}_{\geq 0} = \{ t \in \mathbb{Z} : t \geq 0\}$. Note that any saturation is an ideal in $\mathbb{R}[\bm{z}]$ (\citealp{CLO}; proposition 9, section 4.4). If $W = \mathbb{C}$, the affine variety $\mathbb{V}_{\mathbb{C}}(\mathcal{J}:\mathcal{K}^{\infty})$ coincides with the Zariski closure of the subset $\mathbb{V}_{\mathbb{C}}(\mathcal{J}) \setminus \mathbb{V}_{\mathbb{C}}(\mathcal{K}) \subseteq \mathbb{C}$ (\citealp{CLO}; theorem 10, section 4.4), that is, the $\mathbb{V}_{\mathbb{C}}(\mathcal{J}:\mathcal{K}^{\infty})$ is none other than the Zariski closure of the solution space over $\mathbb{C}$ of 
\begin{align*}
    \left(\eqref{AEGB-1} \mbox{ and } h_1(\bm{z}) \neq 0 \right)
    \mbox{ or }
    \cdots
    \mbox{ or }
    \left(\eqref{AEGB-1} \mbox{ and } h_s(\bm{z}) \neq 0 \right).
\end{align*}
Recall that the Zariski closure of a subset in $W^m = \mathbb{R}^m$ or $\mathbb{C}^m$ is the smallest affine variety containing the subset (see \ref{Fields} for details of Zariski closures).

Now, we provide a proposition that gives theoretical validity to divide Eqs. \eqref{eq:diff1-2} and \eqref{formula-L} into some sub-problems (\citealp{CLO}; theorem 10, section 4.4).

\begin{pro}\label{bunkai}
  \begin{align*}
    \mathbb{V}_W(\mathcal{J}) = \mathbb{V}_W(\mathcal{J} + \mathcal{K}) \cup \mathbb{V}_{W}(\mathcal{J}:\mathcal{K}^{\infty}).
  \end{align*}
\end{pro}
Proposition \ref{bunkai} tells that if we obtain generator sets of the sum $\mathcal{J}+\mathcal{K}$ and the saturation $\mathcal{J}:\mathcal{K}^{\infty}$, and both of their affine varieties are not empty, we get sub-problems of \eqref{AEGB-1}. 

Now, let us turn our attention to the algebraic equation for the maximum likelihood factor analysis. The result of Example \ref{ex1} shows that the affine variety includes many solutions whose $\Psi$ is singular; 6 out of 7 solutions are singular. In fact, as we will see in the numerical example presented in Section \ref{sec:MCS}, the solution space of \eqref{eq:diff1-2} and \eqref{formula-L} includes a surprisingly large number of solutions whose $\Psi$ is singular.  Therefore, we propose to decompose the original problem into the following sub-problems: $\psi_j=0$ and $\psi_j\neq 0$ for $j=1,\dots, p$. Now, we provide a simple example based on Example \ref{ex1} to illustrate how the decomposition can be achieved. A rigorous procedure will be provided in Section \ref{mlsgb-2}.  
\begin{ex}\label{ex1-2}
We illustrate the usefulness of Proposition \ref{bunkai} using the same example as in Example \ref{ex1}. Let us consider the ideal $\mathcal{J} = \langle f_1, \ldots, f_6 \rangle$ in the polynomial ring $\mathbb{R}[\bm{\psi}, \bm{\ell}]$ where $f_1,\ldots,f_6$ are polynomials given in Example \ref{ex1}.  We construct the ideals sequentially as follows: first, by using the ideal $\mathcal{J}$, we construct the following sum and saturation:
  \begin{align*}
    \mathcal{I}_0 = \mathcal{J} + \langle \psi_1 \rangle, ~~ \mathcal{I}_1 = \mathcal{J} : \langle \psi_1 \rangle^{\infty}.
  \end{align*}
  Next, by using the ideals $\mathcal{I}_{0}$ and $\mathcal{I}_{1}$, we construct the following sums and saturations:
  \begin{align*}
    \mathcal{I}_{00} = \mathcal{I}_0 + \langle \psi_2 \rangle, ~~ \mathcal{I}_{01} = \mathcal{I}_0 : \langle \psi_2 \rangle^{\infty},
    \\
    \mathcal{I}_{10} = \mathcal{I}_{1} + \langle \psi_2 \rangle, ~~ \mathcal{I}_{11} = \mathcal{I}_{1} : \langle \psi_2 \rangle^{\infty}.
  \end{align*}
  Last, by using the ideals $\mathcal{I}_{00}$, $\mathcal{I}_{01}$, $\mathcal{I}_{10},$ and $\mathcal{I}_{11}$, we construct the following sums and saturations:
  \begin{align*}
    \mathcal{I}_{000} = \mathcal{I}_{00} + \langle \psi_3 \rangle, ~~
    \mathcal{I}_{001} = \mathcal{I}_{00} : \langle \psi_3 \rangle^{\infty},
    \\
    \mathcal{I}_{010} = \mathcal{I}_{01} + \langle \psi_3 \rangle, ~~
    \mathcal{I}_{011} = \mathcal{I}_{01} : \langle \psi_3 \rangle^{\infty},
    \\
    \mathcal{I}_{100} = \mathcal{I}_{10} + \langle \psi_3 \rangle, ~~
    \mathcal{I}_{101} = \mathcal{I}_{10} : \langle \psi_3 \rangle^{\infty},
\\
    \mathcal{I}_{110} = \mathcal{I}_{11} + \langle \psi_3 \rangle, ~~
    \mathcal{I}_{111} = \mathcal{I}_{11} : \langle \psi_3 \rangle^{\infty}.
\end{align*}  
  It follows from Proposition \ref{bunkai} that
  \begin{align*}
    \mathbb{V}_{\mathbb{R}}(\mathcal{J})
    &= \mathbb{V}_{\mathbb{R}}(\mathcal{I}_0) \cup \mathbb{V}_{\mathbb{R}}(\mathcal{I}_1)
    \\
    &= \left(\mathbb{V}_{\mathbb{R}}(\mathcal{I}_{00}) \cup \mathbb{V}_{\mathbb{R}}(\mathcal{I}_{01})\right) \cup \left(\mathbb{V}_{\mathbb{R}}(\mathcal{I}_{10}) \cup \mathbb{V}_{\mathbb{R}}(\mathcal{I}_{11})\right)
    \\
    &=
    \left(
    \left(\mathbb{V}_{\mathbb{R}}(\mathcal{I}_{000}) \cup \mathbb{V}_{\mathbb{R}}(\mathcal{I}_{001})\right)
    \cup 
    \left(\mathbb{V}_{\mathbb{R}}(\mathcal{I}_{010}) \cup \mathbb{V}_{\mathbb{R}}(\mathcal{I}_{011})\right)
    \right)
    \\
    &\quad
    \cup
    \left(
    \left(\mathbb{V}_{\mathbb{R}}(\mathcal{I}_{100}) \cup \mathbb{V}_{\mathbb{R}}(\mathcal{I}_{101})\right)
    \cup 
    \left(\mathbb{V}_{\mathbb{R}}(\mathcal{I}_{110}) \cup \mathbb{V}_{\mathbb{R}}(\mathcal{I}_{111})\right)
    \right)
    \\
    &=
    \bigcup_{b \in \{ 0,1\}^3} \mathbb{V}_{\mathbb{R}}(\mathcal{I}_{b}).
  \end{align*}
  The affine varieties $\mathbb{V}_{\mathbb{R}}(\mathcal{I}_{000}), \mathbb{V}_{\mathbb{R}}(\mathcal{I}_{100}), \mathbb{V}_{\mathbb{R}}(\mathcal{I}_{010}), \mathbb{V}_{\mathbb{R}}(\mathcal{I}_{001})$ are empty sets. On the other hand,   
  \begin{align*}
    \mathbb{V}_{\mathbb{R}}(\mathcal{I}_{011}) &= \left\{ \left(0, \frac{3}{4}, \frac{8}{9}, \pm 1, \pm\frac{1}{2}, \pm\frac{1}{3}\right) \right\},
    &
    \mathbb{V}_{\mathbb{R}}(\mathcal{I}_{101}) &= \left\{ \left(\frac{3}{4}, 0, \frac{5}{9}, \pm\frac{1}{2}, \pm1, \pm\frac{2}{3}\right)  \right\},
    \\
    \mathbb{V}_{\mathbb{R}}(\mathcal{I}_{110}) &= \left\{ \left(\frac{8}{9}, \frac{5}{9}, 0, \pm\frac{1}{3}, \pm\frac{2}{3}, \pm 1\right) \right\},
    &
    \mathbb{V}_{\mathbb{R}}(\mathcal{I}_{111}) &= \left\{ \left(1, 1, 1, 0, 0, 0\right) \right\}
  \end{align*}
  are non-empty proper subsets of $\mathbb{V}_{\mathbb{R}}(\mathcal{J})$. Thus, the solution space of \eqref{exaegb-1} is divided into the four affine varieties: $\mathbb{V}_{\mathbb{R}}(\mathcal{I}_{011})$, $\mathbb{V}_{\mathbb{R}}(\mathcal{I}_{101})$, $\mathbb{V}_{\mathbb{R}}(\mathcal{I}_{110})$, and  $\mathbb{V}_{\mathbb{R}}(\mathcal{I}_{111})$.
\end{ex}

\subsubsection{Gr\"obner bases}
Example \ref{ex1-2} showed that the solution space can be successfully decomposed by constructing sums and saturations. To obtain some simple sub-problems of the algebraic equations \eqref{eq:diff1-2} and \eqref{formula-L} for maximum likelihood factor analysis, in addition, we need to compute generator sets of a sum $\mathcal{J} + \mathcal{K}$ and a saturation $\mathcal{J}:\mathcal{K}^{\infty}$. The sum $\mathcal{J} + \mathcal{K}$ is generated by $f_1, \ldots, f_r, h_1, \ldots, h_s$, that is, $\mathcal{J} + \mathcal{K} = \langle f_1,\ldots, f_r, h_1,\ldots,h_s \rangle$. A generator set of the saturation $\mathcal{J}:\mathcal{K}^{\infty}$ is obtained by using a Gr\"obner basis. The definition of the Gr\"obner basis requires a monomial ordering. Therefore, we define monomial orderings before giving the definition of Gr\"obner bases. 

Let us denote $\bm{z}^{\bm{a}} = z_1^{a_1} \cdots z_m^{a_m}$ for $\bm{a}=(a_1,\ldots,a_m) \in \mathbb{Z}_{\geq 0}^m$, which is called a monomial, where $\mathbb{Z}_{\geq 0}^m = \{ (\alpha_1,\ldots,\alpha_m) \in \mathbb{Z}^m : \alpha_1, \ldots, \alpha_m \geq 0\}$. Let $M(\bm{z}) = \{ \bm{z}^{\bm{a}} : \bm{a} \in \mathbb{Z}_{\geq 0}^m \}$.
\begin{dfn}
  A {\em monomial ordering} $\succ$ on $\mathbb{R}[\bm{z}]$ is a relation on the set $M(\bm{z})$ satisfying that  
  \begin{enumerate}
  \item $\succ$ is a total ordering on $M(\bm{z})$,  
  \item $\bm{z}^{\bm{a}} \succ \bm{z}^{\bm{b}} \Longrightarrow \bm{z}^{\bm{a}} \bm{z}^{\bm{c}} \succ \bm{z}^{\bm{b}} \bm{z}^{\bm{c}}$ for any $\bm{c} \in \mathbb{Z}_{\geq 0}^m$, and 
  \item any nonempty subset of $M(\bm{z})$ has a smallest element under $\succ$.
  \end{enumerate}
  Fixing a monomial ordering on $\mathbb{R}[\bm{z}]$ and given a polynomial $f = \sum_{\bm{a} \in \mathbb{Z}_{\geq 0}^m} c_{\bm{a}} \bm{z}^{\bm{a}}$ for $c_{\bm{a}} \in \mathbb{R}$, we call the monomial $\mathrm{LM}(f) = \max(\bm{z}^{\bm{a}} : c_{\bm{a}} \neq 0)$ the leading monomial of $f$.
\end{dfn}

We give two examples of monomial orderings. The first one is called a {\em lexicographic order} (or the {\em lex order} for short). The second one is called a {\em graded reverse lex order} (or the {\em grevlex order} for short). As shown in the following example, the leading monomial of a given polynomial is determined depending on a fixed monomial ordering. 
\begin{ex}[Lexicographic Order]
  Let $\bm{a}, \bm{b} \in \mathbb{Z}_{\geq 0}^m$. We say $\bm{z}^{\bm{a}} \succ_{\mathrm{lex}} \bm{z}^{\bm{b}}$ if the leftmost nonzero entry of the vector difference $\bm{a} - \bm{b}$ is positive.
  For example, when we fix the lex order, $\mathrm{LM}(3 z_1^2 + 2 z_1 z_2^3) = z_1^2$. In addition,
  \begin{align*}
    \mathrm{LM}(f_4) = \psi_1, ~~ \mathrm{LM}(f_5) = \psi_2, ~~ \mathrm{LM}(f_6) = \psi_3,
  \end{align*}
  where $f_4,f_5,f_6 \in \mathbb{R}[\bm{\psi},\bm{\ell}]$ are the polynomials presented in the Example \ref{ex1}.
\end{ex}

\begin{ex}[Graded Reverse Lex Order]
  Let $\bm{a} = (a_1,\ldots, a_m), \bm{b} = (b_1,\ldots, b_m) \in \mathbb{Z}_{\geq 0}^m$. We say $\bm{z}^{\bm{a}} \succ_{\mathrm{grevlex}} \bm{z}^{\bm{b}}$ if $|\bm{a}| = \sum_{i=1}^{m} a_i > |\bm{b}| = \sum_{i=1}^{m} b_i$ or the rightmost nonzero entry of the vector difference $\bm{a} - \bm{b}$ is negative. For example, when we fix the grevlex order, $\mathrm{LM}(3 z_1^2 + 2 z_1 z_2^3) = z_1 z_2^3$. In addition,
  \begin{align*}
\mathrm{LM}(f_4) = \ell_{11}^3, ~~ \mathrm{LM}(f_5) = \ell_{11}^2 \ell_{21}, ~~ \mathrm{LM}(f_6) = \ell_{11}^2 \ell_{31},
  \end{align*}
  where $f_4,f_5,f_6 \in \mathbb{R}[\bm{\psi},\bm{\ell}]$ are the polynomials presented in the Example \ref{ex1}.
\end{ex}

Now we define Gr\"obner bases.
\begin{dfn}
  Fix a monomial ordering on $\mathbb{R}[\bm{z}]$. A finite subset $G = \{ g_1, \ldots, g_t \}$ of an ideal $\mathcal{J} \subseteq \mathbb{R}[\bm{z}]$ is said to be a {\em Gr\"obner basis} of $\mathcal{J}$ if
  \begin{align*}
    \langle \mathrm{LM}(g_1), \ldots, \mathrm{LM}(g_t) \rangle
    =
    \langle \mathrm{ \{ LM}(f) : f \in \mathcal{J}  \}  \rangle
  \end{align*}
\end{dfn}
Note that every polynomial ideal has a Gr\"obner basis (\citealp{CLO}; corollary 6, section 2.5), and that any Gr\"obner basis of $\mathcal{J}$ generates the ideal $\mathcal{J}$.

For the case $\mathcal{K} = \langle h \rangle \subseteq \mathbb{R}[\bm{z}]$, a generator set of the saturation $\mathcal{J}:\mathcal{K}^{\infty}$ is obtained by using a Gr\"obner basis as follows (\citealp{CLO}; theorem 14 (ii), section 4.4):
\begin{enumerate}
\item let $\tilde{\mathcal{J}} = \langle f_1, \ldots, f_r, 1 - y h \rangle \subseteq \mathbb{R}[y, \bm{z}]$,
\item compute a Gr\"obner basis $\tilde{G}$ of $\tilde{\mathcal{J}}$ with respect to the lex order,
\item then $G = \tilde{G} \cap \mathbb{R}[\bm{z}]$ is a Gr\"obner basis of the saturation $\mathcal{J}:\mathcal{K}^{\infty}$.
\end{enumerate}

\begin{ex}\label{ex1-3}
Let us reconsider Example \ref{ex1}, which is related to the algebraic equations \eqref{eq:diff1-2} and \eqref{formula-L}. Gr\"obner bases of the sums and saturations described above give us simplified sub-problems.  
  For $b \in \{0,1\}^3$, a Gr\"obner basis $G_b$ of $\mathcal{I}_{b}$ with respect to the lex order is given by 
  \begin{align*}
    G_{000} &= \left\{ 1 \right\},~~G_{001} = \left\{ 1 \right\},~~G_{010} = \left\{ 1 \right\},~~G_{100} = \left\{ 1 \right\},
    \\
    G_{011} &= \left\{ \psi_1, \psi_2-\frac{3}{4}, \psi_3-\frac{8}{9}, \ell_{11}-3 \ell_{31}, \ell_{21} - \frac{3}{2} \ell_{31}, \ell_{31}^2-1/9 \right\},
    \\
    G_{101} &= \left\{ \psi_1-\frac{3}{4}, \psi_2, \psi_3-\frac{5}{9}, \ell_{11}-\frac{3}{4} \ell_{31}, \ell_{21} - \frac{3}{2} \ell_{31}, \ell_{31}^2-\frac{4}{9}  \right\},
    \\
    G_{110} &= \left\{ \psi_1-\frac{8}{9}, \psi_2-\frac{5}{9}, \psi_3, \ell_{11}-\frac{1}{3} \ell_{31}, \ell_{21} - \frac{2}{3} \ell_{31}, \ell_{31}^2-1  \right\},
    \\
    G_{111} &= 
    \left\{
    \psi_1+\frac{81}{25} \ell_{31}^2-1, \psi_2-1, \psi_3+\ell_{31}^2-1, \ell_{11}+\frac{9}{5} \ell_{31}, \ell_{21}, \ell_{31}^3 + \frac{5}{27} \ell_{31}
    \right\}.
  \end{align*}
  So the Gr\"obner bases $G_{011}, G_{101}, G_{110}, G_{111}$ divide \eqref{exaegb-1} into the following simple sub-problems:
  \begin{align*}
    &
    \left\{\begin{array}{l}
      \psi_1=0,
      \\
      \psi_2=\frac{3}{4},
      \\
      \psi_3=\frac{8}{9},
      \\
      \ell_{11}=3 \ell_{31},
      \\
      \ell_{21}=\frac{3}{2} \ell_{31},
      \\
      0 = \ell_{31}^2-1/9,
    \end{array}\right.
    &&
    \left\{\begin{array}{l}
      \psi_1=\frac{3}{4},
      \\
      \psi_2=0,
      \\
      \psi_3=\frac{5}{9},
      \\
      \ell_{11}=\frac{3}{4} \ell_{31},
      \\
      \ell_{21} = \frac{3}{2} \ell_{31},
      \\
      0 = \ell_{31}^2-\frac{4}{9},
    \end{array}\right.
    \\
    &
    \left\{\begin{array}{l}
      \psi_1=\frac{8}{9},
      \\
      \psi_2=\frac{5}{9},
      \\
      \psi_3=0,
      \\
      \ell_{11}=\frac{1}{3} \ell_{31},
      \\
      \ell_{21} = \frac{2}{3} \ell_{31},
      \\
      0=\ell_{31}^2-1,
    \end{array}\right.
    &&
    \left\{\begin{array}{l}
      \psi_1=1-\frac{81}{25} \ell_{31}^2,
      \\
      \psi_2=1,
      \\
      \psi_3=1-\ell_{31}^2,
      \\
      \ell_{11}=-\frac{9}{5} \ell_{31},
      \\
      \ell_{21}=0,
      \\
      0 = \ell_{31}^3 + \frac{5}{27} \ell_{31}.
    \end{array}\right.
  \end{align*}
\end{ex}

Gr\"obner bases have a good property that plays an important role in solving an algebraic equation of the form \eqref{AEGB-1}. The following proposition shows that Gr\"obner bases transform \eqref{AEGB-1} into easily solvable problems, without being subject to chance. The following proposition is called the finiteness theorem (\citealp{CLO}; theorem 6, section 5.3).

\begin{pro}\label{finiteness-thm}
  Fix a monomial ordering on $\mathbb{R}[\bm{z}]$ and let $G$ be a Gr\"obner basis of an ideal $\mathcal{J} \subseteq \mathbb{R}[\bm{z}]$.
  Consider the following four statements:
  \begin{enumerate}
  \item\label{finiteness-thm-1} for each $i = 1, \ldots, m$, there exists $t_i \geq 0$ such that $z_i^{t_i} \in \langle \mathrm{LM}(f) : f \in \mathcal{J} \rangle$,
  \item\label{finiteness-thm-2} for each $i = 1, \ldots, m$, there exists $t_i \geq 0$ and $g \in G$ such that $\mathrm{LM}(g) = z_i^{t_i}$,
  \item\label{finiteness-thm-3} for each $i = 1, \ldots, m$, there exists $t_i \geq 0$ and $g \in G$ such that $\mathrm{LM}(g) = z_i^{t_i}$ and $g \in \mathbb{R}[z_i, \ldots, z_m]$, if we fix the lex order,
  \item\label{finiteness-thm-4} the affine variety $\mathbb{V}_{W}(\mathcal{J})$ is a finite set.
  \end{enumerate}
  Then the statements \ref{finiteness-thm-1}-\ref{finiteness-thm-3} are equivalent and they all imply the statement \ref{finiteness-thm-4}. An ideal $\mathcal{J}$ satisfying the statement \ref{finiteness-thm-1}, \ref{finiteness-thm-2}, or \ref{finiteness-thm-3} is called a {\em zero-dimensional ideal}. Otherwise it is called a {\em non-zero dimensional ideal}. Furthermore, if $W = \mathbb{C}$, then the statements \ref{finiteness-thm-1}-\ref{finiteness-thm-4} are all equivalent. 
\end{pro}

The above proposition shows that if we can obtain a Gröbner basis, we can determine whether the affine variety is a finite or infinite set, and if the affine variety is a finite set, the Gröbner basis with respect to the lex order gives a triangulated representation (statement \ref{finiteness-thm-3} of Proposition \ref{finiteness-thm}).  This representation produces a much simpler algebraic equation than the original one; for example, when the affine variety is a finite linear set, the Gröbner basis provides a triangular matrix. In general, the computation of a Gr\"obner basis of the grevlex order is more efficient than that of the lex order. Also, if $\mathcal{J}$ is zero-dimensional, we have an efficient algorithm which converts a Gr\"oebner basis from one monomial orderings to another, which is called the FGLM algorithm \citep{FGLM}. Hence, in our algebraic approach to maximum likelihood factor analysis, we first compute Gr\"obner bases with respect to the grevlex order, and then convert them to the lex order.

Primary ideal decompositions can {\it more completely} give us simplified sub-problems of algebraic equations. However, unfortunately, the primary ideal decompositions require {\it heavy} computational loads in general. Thus, we employ Proposition \ref{bunkai}  instead of primary ideal decompositions.

\subsection{An algorithm to find all solutions}\label{mlsgb-2}

We provide an algorithm to find all maximum likelihood solutions to Eq. \eqref{eq:diff0}. To get the solution, we solve the following equation:
\begin{align}
\{S-LL^{\top}-\mathrm{diag}(S-LL^{\top})\}S^{-1}L =0,\label{Ldiff3}
\end{align}
which is derived by substituting \eqref{eq:diff1-2} into \eqref{formula-L}. Note that, according to Theorem \ref{we-want-solve}, the algebraic equation \eqref{Ldiff3} is a necessary condition for the likelihood equation \eqref{eq:diff0}.  Since the solution space of \eqref{Ldiff3} contains that of \eqref{eq:diff0}, we can find all solutions to \eqref{eq:diff0} by solving \eqref{Ldiff3}.  Eq. \eqref{Ldiff3} is much easier to solve than Eq. \eqref{eq:diff0}; this is because Eq. \eqref{eq:diff0} includes an inverse $\Sigma^{-1}$ that consists of rational functions, while   Eq. \eqref{Ldiff3} does not.

As described in Example \ref{ex1-2} in the section \ref{mlsgb-1}, \eqref{Ldiff3} can have a large number of, and sometimes infinite, solutions when $\Psi$ is singular (see \ref{detail-implimentation} for details).  Therefore, as in Example \ref{ex1-2}, we use Proposition \ref{bunkai} to get some sub-problems of Eq. \eqref{Ldiff3}. 

Considering a rotational indeterminacy, we suppose that the upper triangular elements of $L$ are fixed by zero, that is, 
\begin{align*}
  L
  =
  \begin{pmatrix}
    \ell_{11}                                           \\
    \ell_{21} & \ell_{22} &        & \text{\huge{0}}    \\
    \vdots    & \vdots    & \ddots                      \\
    \vdots    & \vdots    &        & \ddots             \\
    \ell_{k1} & \ell_{k2} & \cdots & \cdots & \ell_{kk} \\
    \vdots    & \vdots    & \vdots & \vdots & \vdots    \\
    \ell_{p1} & \ell_{p2} & \cdots & \cdots & \ell_{pk}     
\end{pmatrix}.
\end{align*}
Let us consider the following ideal in $\mathbb{R}[\bm{\ell}] = \mathbb{R}[\ell_{ij} : 1 \leq i \leq p, 1 \leq j \leq \min(i,k)]$:
\begin{align*}
  \mathcal{J} = \langle (\{S-LL^{\top}-\mathrm{diag}(S-LL^{\top})\}S^{-1}L)_{ij} : 1 \leq i \leq p, 1 \leq j \leq \min(i,k) \rangle.
\end{align*}
We have to choose the ideals of the form $\mathcal{K} = \langle h \rangle \subseteq \mathbb{R}[\bm{\ell}]$ well to get sub-problems of \eqref{Ldiff3}.

Recall that \eqref{eq:diff1-2} is equivalent to the lower part of \eqref{eq:diff0} when $\Psi$ is not singular, and is a necessary condition of \eqref{eq:diff0} when $\Psi$ is singular as shown in Theorem \ref{we-want-solve}. Focusing on \eqref{eq:diff1-2}, we consider $\psi_1, \ldots, \psi_p$ as polynomials defined by
\begin{align*}
  \psi_{i} = s_{ii} - \sum_{j=1}^{\min(i,k)} \ell_{ij}^2 \in \mathbb{R}[\bm{\ell}] ~~ (i=1,\ldots,p).
\end{align*}
For each $\mathcal{K}_{i} = \langle \psi_i \rangle \subseteq \mathbb{R}[\bm{\ell}]$, we construct sums and saturations as follows.

First, we construct the sum $\mathcal{I}_0 = \mathcal{J} + \mathcal{K}_{1}$ and the saturation $\mathcal{I}_1 = \mathcal{J}:\mathcal{K}_{1}^{\infty}$. Proposition \ref{bunkai} implies that
\begin{align*}
  \mathbb{V}_{W}(\mathcal{J}) = \mathbb{V}_{W}(\mathcal{I}_0) \cup \mathbb{V}_{W}(\mathcal{I}_1).
\end{align*}
Second, we construct the follwoing sums and saturations:
\begin{align*}
  \mathcal{I}_{00} = \mathcal{I}_{0} + \mathcal{K}_{2}, ~~ \mathcal{I}_{01} = \mathcal{I}_0:\mathcal{K}_{2}^{\infty}, ~~ \mathcal{I}_{10} = \mathcal{I}_{1} + \mathcal{K}_{2}, ~~ \mathcal{I}_{11} = \mathcal{I}_1:\mathcal{K}_{2}^{\infty}.
\end{align*}
Since Proposition \ref{bunkai} implies $\mathbb{V}_{W}(\mathcal{I}_0) = \mathbb{V}_{W}(\mathcal{I}_{00}) \cup \mathbb{V}_{W}(\mathcal{I}_{01})$ and $\mathbb{V}_{W}(\mathcal{I}_1) = \mathbb{V}_{W}(\mathcal{I}_{10}) \cup \mathbb{V}_{W}(\mathcal{I}_{11})$, 
\begin{align*}
  \mathbb{V}_{W}(\mathcal{J}) = \mathbb{V}_{W}(\mathcal{I}_{00}) \cup \mathbb{V}_{W}(\mathcal{I}_{01}) \cup \mathbb{V}_{W}(\mathcal{I}_{10}) \cup \mathbb{V}_{W}(\mathcal{I}_{11})
\end{align*}
holds. In the same way, for each $b \in \{ 0, 1\}^{i}$, $2 \leq i \leq p - 1$, we sequentially construct
\begin{align*}
  \mathcal{I}_{b0} = \mathcal{I}_{b} + \mathcal{K}_{i+1}, ~~ \mathcal{I}_{b1} = \mathcal{I}_b:\mathcal{K}_{i+1}^{\infty}.
\end{align*}
By Proposition \ref{bunkai}, we have $\mathbb{V}_{W}(\mathcal{J}) = \bigcup_{b \in \{ 0, 1\}^{p}} \mathbb{V}_{W}(\mathcal{I}_{b})$. Finally, we compute all solutions to \eqref{Ldiff3} by using the following steps for each $b \in \{ 0, 1\}^{p}$:
\begin{enumerate}
\item compute a Gr\"obner basis $\tilde{G}_b = \{ \tilde{g}_1, \ldots, \tilde{g}_s\}$ of $\mathcal{I}_{b}$ with respect to the grevlex order,
\item if $\mathcal{I}_{b}$ is zero-dimensional, compute a Gr\"obner basis $G_b =  \{ g_1, \ldots, g_{r}\}$ of $\mathcal{I}_b$ with respect to the lex order by converting $\tilde{G}_b$ to the lex order, and then solve the equation $g_1 = \cdots = g_r = 0$,  
\item else, compute sample points of each connected component defined by $\tilde{g}_1 = \cdots = \tilde{g}_s = 0$ by using cylindrical decomposition.
\end{enumerate}
The above algorithm is detailed in Algorithm \ref{alg:fukafuka}.  
\begin{rem}
  In practice, the algebraic equation \eqref{Ldiff3} is divided into subproblems by also using polynomials other than $\psi_1, \ldots, \psi_p \in \mathbb{R}[\bm{\ell}]$ because the algebraic equation \eqref{Ldiff3} has a large number of  solutions in general. In addition, whether $\mathcal{I}_{b}$ is zero-dimensional or not, it can be determined by the statement \ref{finiteness-thm-2} of Proposition \ref{finiteness-thm}. Since we divide \eqref{Ldiff3} into a large number of sub-problems, even if $\mathcal{I}_b$ is non zero-dimensional, $\tilde{G}_b$ has a simple representation. If $\mathcal{I}_{b}$ is zero-dimensional, the Gr\"obner basis $G_b$ gives us an easily solvable sub-problem, as described in Example \ref{ex1-3} in Section \ref{mlsgb-2} and the statement \ref{finiteness-thm-3} of Proposition \ref{finiteness-thm}. For details of the decomposition other than $\psi_1, \ldots, \psi_p \in \mathbb{R}[\bm{\ell}]$, please refer to \ref{detail-implimentation}.
\end{rem}
\begin{rem}
  The ideals $\mathcal{I}_b$ may have multiple roots. Therefore, in practice, Step 1 computes a Gr\"obner basis of the radical ideals for $\mathcal{I}_{b}$. So our algebraic approach can treat very simple sub-problems, and does not generate multiple solutions.
\end{rem}
\begin{algorithm}[t]
\caption{An algorithm to find all solutions}
    \label{alg:fukafuka}
    \begin{algorithmic}[1]    
    \REQUIRE Sample variance-covariance matrix $S = (s_{ij})$ 
    \ENSURE All solutions to the algebraic equations \eqref{eq:diff1-2} and \eqref{formula-L} ({\tt All.Sol.})
    \STATE $\mathcal{J}$ $\leftarrow$ the ideal $\langle (\{S-LL^{\top}-\mathrm{diag}(S-LL^{\top})\}S^{-1}L)_{ij} : 1 \leq i \leq p, 1 \leq j \leq  \min(i,k) \rangle$
    \STATE $h_i$ $\leftarrow$ the polynomial $s_{ii} - \sum_{j=1}^{k} \ell_{ij}^2$ for $1 \leq i \leq p$
    \STATE $\mathcal{K}_{i}$ $\leftarrow$ the ideal $\langle h_i \rangle$ for $1 \leq i \leq p$
    \STATE $\mathcal{I}_0$ $\leftarrow$ the sum $\mathcal{J} + \mathcal{K}_{1}$
    \STATE $\mathcal{I}_1$ $\leftarrow$ the saturation $\mathcal{J}:\mathcal{K}_{1}^{\infty}$
    \FOR{$1 \leq i \leq p  - 1$}
    	\FOR{$b \in \{ 0, 1\}^{i-1}$}
    	\STATE $\mathcal{I}_{b0}$ $\leftarrow$ the sum $\mathcal{I}_{b} + \mathcal{K}_{i+1}$
        \STATE $\mathcal{I}_{b1}$ $\leftarrow$ the saturation $\mathcal{I}_b:\mathcal{K}_{i+1}^{\infty}$
        \ENDFOR
    \ENDFOR
    \FOR{$b \in \{ 0, 1\}^{p}$}
        \STATE $\tilde{G}_b = \{ \tilde{g}_1, \ldots, \tilde{g}_s\}$ $\leftarrow$ a Gr\"obner basis of $\mathcal{I}_{b}$ with respect to the grevlex order
        \IF{$\mathcal{I}_b$ is zero-dimensional}
            \STATE $G_b = \{ g_1, \ldots, g_r\}$ $\leftarrow$ a Gr\"obner basis of $\mathcal{I}_b$ with respect to the lex order by converting $\tilde{G}_b$
            \STATE {\tt All.Sol.} $\leftarrow$ {\tt All.Sol.} $\cup$ $\{ (L, \mathrm{diag}(S - LL^\top)) \in \mathbb{R}^{p \times k} \times \mathbb{R}^{p \times p} : g_1 = \cdots = g_r = 0\}$
        \ELSE
            \STATE {\tt All.Sol.} $\leftarrow$ {\tt All.Sol.} $\cup$ $\left\{\begin{array}{r}(L, \mathrm{diag}(S - LL^\top)) \in \mathbb{R}^{p \times k} \times \mathbb{R}^{p \times p} : \\ \bm{\ell} \mbox{ is a sample points of}\\ \mbox{each connected component}\\ \mbox{defined by } \tilde{g}_1 = \cdots = \tilde{g}_s = 0\end{array}\right\}$
        \ENDIF
    \ENDFOR
    \end{algorithmic}
\end{algorithm}

\subsection{An algorithm to identify a solution pattern}\label{mlsgb-3}
 We provide an algorithm to categorize the maximum likelihood solution into three patterns: ``Proper solution", ``Improper solution", and ``No solutions". The algorithm is detailed in Algorithm \ref{alg:detectimproper_algsol}. First, we find all solutions of Eq. \eqref{eq:diff0} with Algorithm \ref{alg:fukafuka}. Then, we get a solution that minimizes the discrepancy function in  \eqref{q}.  When the Hessian matrix at the solution is positive definite, there exists a maximum likelihood solution that minimizes \eqref{q}. In this case, the solution is categorized as either proper or improper according to the sign of the unique variances. When all unique variances are positive, the solution is referred to as ``proper solution". Meanwhile, when at least one of the unique variances is not positive, the solution is referred to as a ``improper solution".  When the Hessian matrix at the solution is not positive definite, there are no optimal maximum likelihood solutions. With Algorithm \ref{alg:detectimproper_algsol}, we can investigate the tendency of the solution pattern through Monte Carlo simulations, which will be described in the following section. 

\begin{algorithm}[t]
\caption{An algorithm to identify the solution pattern}
    \label{alg:detectimproper_algsol}
    \begin{algorithmic}[1]    
    \REQUIRE All solutions of \eqref{Ldiff3}, say $(\hat{L}_{t}, \hat{\Psi}_{t})_{t=1,\dots,T}.$ 
    \ENSURE Solution pattern ({\tt Sol.Pat.})
    \STATE Select a solution that minimizes the discrepancy function in Eq. \eqref{q}; that is, obtain ($\hat{L}_{t^*}, \hat{\Psi}_{t^*}$), where $\displaystyle {t^*} = \argmin_{t \in \{1,\dots,T\}} \ell(\hat{L}_t, \hat{\Psi}_t)$.  
    \IF{the solution ($\hat{L}_{t^*}, \hat{\Psi}_{t^*}$) satisfies the condition \eqref{eq:diff0} and its observed Fisher information matrix is positive-definite}
    	\IF{the minimum value of diagonal element of $\hat{\Psi}_{t^*}$ is positive}
    	\STATE {\tt Sol.Pat.} $\leftarrow$ ``Proper solution."
    	\ELSE
    	\STATE {\tt Sol.Pat.} $\leftarrow$ ``Improper solution."    	\ENDIF
\ELSE
    	\STATE {\tt Sol.Pat.} $\leftarrow$ ``No solutions."
    \ENDIF
    
    \end{algorithmic}
\end{algorithm}

 \section{Monte Carlo simulation}\label{sec:MCS}
 We conduct Monte Carlo simulations to investigate the characteristics of solution patterns of maximum likelihood estimate in various situations.  First, we consider the following three simulation models:
 \begin{align*}
  {{\rm S1}: \ \ }
   L = 
    \begin{pmatrix}
      0.9 & 0
      \\
      0.8 & 0
      \\
      0.7 & 0
      \\
      0  & 0.8
      \\
      0  & 0.7
    \end{pmatrix},
    \quad
  {{\rm S2}: \ \ }
    L = 
    \begin{pmatrix}
      0.5 & 0
      \\
      0.5 & 0
      \\
      0.5 & 0.5
      \\
      0  & 0.5
      \\
      0  & 0.5
    \end{pmatrix},
    \quad
  {{\rm S3}: \ \ }
    L = 
    \begin{pmatrix}
      0.9 & 0
      \\
      0.8 & 0
      \\
      0.6 & 0.7
      \\
      0  & 0.8
      \\
      0  & 0.9
    \end{pmatrix},
\end{align*}
 and $\Psi={\rm diag}(I-LL^\top)$. On S1, the second column of $L$ has only two nonzero elements, suggesting that a necessary condition for identification is not satisfied (\citealp{anderson1956statistical}; theorem 5.6).  In this case, the estimate can be unstable due to the identifiability problem \citep{Driel.1978}, and improper solutions are often obtained.  S2 has small communalities and large unique variances, and then improper solutions sometimes occur due to the sampling fluctuation \citep{Driel.1978}. In S3, communalities are large and a necessary condition on identification is satisfied; thus, it is expected that numerical estimate is stable and improper solutions are unlikely to be obtained.   
 
For all simulation models, we generate $N=100$ observations from $N(\bm{0}, LL^\top+\Psi)$, and compute the sample covariance matrix.  The elements of the sample covariance matrix are rounded to one decimal place because numbers with many decimal places require a significant amount of computer memory with computational algebra.  For each simulation model, we run 10 times and investigate the characteristics of exact solutions obtained by computational algebra. We note that such a small number of replications of the Monte Carlo simulation is due to the heavy computational loads of computational algebra.  For implementation, we used the highest-performance computer that we had (Intel Xeon Gold 6134 processors, 3.20GHz, 32 CPUs, 256 GB memory, Ubuntu 18.04.2 LTS), and made efforts to optimize the algorithm of computational algebra as much as possible, as detailed in Appendix B. However, the computational timing varied from 2 weeks to 4 weeks to get the result for a single dataset. For comparison, we also perform \factanal function in \R.

\begin{table}[t]
  \caption{Solution pattern obtained by Algorithm \ref{alg:detectimproper_algsol}: proper solution (P), improper solution (I), and no solution (NA).}
  \label{tab:pattern}
\begin{center}
    \small
  \begin{tabular}{lcccccccccc}
\hline
    Sim\#    & 1     & 2    & 3    & 4     & 5    & 6    & 7    &8    & 9     &10   \\ \hline
    S1    & P     & NA    & NA    & NA    & NA    & I     & I     & P     & P     & P \\
    S2    & P     & P     & I     & NA    & P     & I     & I     & P     & P     & P \\
    S3    & P     & P     & P     & P     & P     & P     & P     & P     & P     & P \\
\hline
    \end{tabular}
\end{center}
\end{table}

\subsection{Investigation of maximum likelihood solutions}
We investigate the maximum likelihood solutions obtained by computational algebra and {\tt factanal}.  Table \ref{tab:pattern} shows the solution pattern of the maximum likelihood estimate computed using Algorithm \ref{alg:detectimproper_algsol}. For S3, we obtain proper solutions in all 10 simulations, implying that the results for S3 are stabler than S1 and S2.  Meanwhile, several improper solutions are obtained, or no solutions are found for S1 and S2.  In particular, S1 tends to produce ``no solutions" more frequently than S2.

  With \factanal function, however, it would be difficult to identify the solution pattern as in table \ref{tab:pattern} because the \factanal provides a solution under a restriction that all unique variances are smaller than some threshold, such as 0.005. Moreover, the observed Fisher information matrices of the maximum likelihood solution computed with \factanal function are positive definite for all simulated datasets; therefore, it is remarkably challenging to distinguish ``improper solution" and ``no solutions" using \factanal function.

To study the solution pattern in more detail, we investigate the number of solutions obtained by computational algebra, as shown in Table \ref{tab:numsol}. When counting the number of solutions, the identifiability with respect to the sign of column vectors in a loading matrix is addressed.  
\begin{table}[!t]
  \caption{  \small
The number of solutions obtained by computational algebra. When counting the number of solutions, the identifiability with respect to the sign of column vectors in a loading matrix is addressed.}
  \label{tab:numsol}
\begin{center}
  \small
    \begin{tabular}{lrrrrrrrrrl}
\hline
    Sim\#    & 1     & 2    & 3    & 4     & 5    & 6    & 7    &8    & 9     &10   \\ \hline
 \multicolumn{11}{p{9cm}}{(a) All solutions to \eqref{eq:diff1-2} and \eqref{formula-L} whose unique variances are all positive.}   \\ \hline
    S1    & 3     & 1     & 2     & 3     & 3     & 3     & 0     & 3     & 1     & 1 \\
    S2    & 8     & 2     & 3     & 5     & 3     & 0     & 2     & 2     & 3     & 4 \\
    S3    & 1     & 1     & 3     & 2     & 1     & 1     & 2     & 1     & 2     & 1 \\ 
 \hline\hline
 \multicolumn{11}{p{9cm}}{(b) Solutions of (a) whose observed Fisher information matrix is positive definite.}    \\ \hline
    S1    & 1     & 0     & 0     & 0     & 0     & 0     & 0     & 1     & 1     & 1 \\
    S2    & 1     & 1     & 0     & 0     & 1     & 0     & 0     & 1     & 1     & 1 \\
    S3    & 1     & 1     & 1     & 1     & 1     & 1     & 1     & 1     & 1     & 1 \\
 \hline\hline
 \multicolumn{11}{p{9cm}}{(c) All solutions to \eqref{eq:diff1-2} and \eqref{formula-L} where  all unique variances are non-zero and some of them are negative.}    \\ \hline
    S1    & 2     & 7     & 4     & 2     & 3     & 0     & 10    & 0     & 3     & 0 \\
    S2    & 4     & 3     & 2     & 5     & 3     & 4     & 2     & 2     & 3     & 2 \\
    S3    & 0     & 6     & 1     & 1     & 2     & 4     & 3     & 7     & 8     & 5 \\
 \hline\hline
 \multicolumn{11}{p{9cm}}{(d) Solutions of (c) whose observed Fisher information matrix is positive definite.}    \\ \hline
    S1    & 0     & 0     & 0     & 0     & 0     & 0     & 1     & 0     & 0     & 0 \\
    S2    & 0     & 0     & 2     & 0     & 0     & 1     & 1     & 0     & 0     & 0 \\
    S3    & 0     & 0     & 0     & 0     & 0     & 0     & 0     & 0     & 0     & 0 \\
 \hline\hline
 \multicolumn{11}{p{9cm}}{(e) All solutions to \eqref{eq:diff1-2} and \eqref{formula-L} where some of the unique variances are zero.   }    \\ \hline
    S1    & 31    & 28    & 32    & 30    & 30    & 22    & 32    & 28    & 24    & 26 \\
    S2    & 31    & 22    & 26    & 28    & 27    & 23    & 24    & 27    & 26    & 25 \\
    S3    & 22    & 22    & 23    & 24    & 23    & 24    & 27    & 23    & 28    & 25 \\
 \hline\hline
 \multicolumn{11}{p{9cm}}{(f) Solutions of (e) whose observed Fisher information matrix is positive definite.}    \\ \hline
    S1    & 1     & 5     & 4     & 1     & 3     & 5     & 1     & 2     & 2     & 0 \\
    S2    & 2     & 3     & 6     & 2     & 0     & 3     & 1     & 1     & 1     & 3 \\
    S3    & 1     & 4     & 3     & 1     & 1     & 1     & 2     & 2     & 2     & 3 \\
 \hline\hline
 \multicolumn{11}{p{9cm}}{(g) Solutions of (e) that satisfy \eqref{eq:diff0}.}    \\ \hline
S1 & 0 & 0 & 0 & 0 & 0 & 1 & 1 & 0 & 0 & 4 \\
S2 & 0 & 0 & 0 & 0 & 0 & 0 & 1 & 0 & 0 & 0 \\
S3 & 4 & 0 & 0 & 0 & 0 & 0 & 1 & 0 & 0 & 0 \\
 \hline\hline
 \multicolumn{11}{p{9cm}}{(h) Solutions of (f) that satisfy \eqref{eq:diff0}.}    \\ \hline
S1 & 0 & 0 & 0 & 0 & 0 & 1 & 0 & 0 & 0 & 0 \\
S2 & 0 & 0 & 0 & 0 & 0 & 0 & 0 & 0 & 0 & 0 \\
S3 & 0 & 0 & 0 & 0 & 0 & 0 & 0 & 0 & 0 & 0 \\
 \hline
\end{tabular}
\end{center}
\end{table}
For all cases, the number of solutions that satisfy the conditions \eqref{eq:diff1-2} and \eqref{formula-L} are relatively large, but most of them do not have a positive-definite observed Fisher information matrix.   Thus, it would be essential to check if the observed Fisher information matrix is positive or not. Another noteworthy point is that although we get a number of solutions with zero unique variances using \eqref{eq:diff1-2} and \eqref{formula-L}, about 98.5 \% of them are not the solution of \eqref{eq:diff0} (see the results of (e) and (f)), implying that most of the solutions to \eqref{eq:diff1-2} and \eqref{formula-L} with zero unique variances are inappropriate. As previously mentioned, the conditions \eqref{eq:diff1-2} and \eqref{formula-L} are equivalent to a condition \eqref{eq:diff0} under the assumption that $\Psi$ is not singular, while \eqref{eq:diff1-2} and \eqref{formula-L} are necessary conditions for \eqref{eq:diff0} when $\Psi$ is singular; thus, the solution to \eqref{eq:diff1-2} and \eqref{formula-L} whose unique variances are singular cannot always satisfy \eqref{eq:diff0}.  
%Furthermore, some of the solutions of \eqref{eq:diff1-2} and \eqref{formula-L} that are not the solutions of \eqref{eq:diff0} have a positive-definite observed Fisher information matrix (see the result (f)). Therefore, the conventional numerical algorithm based on \eqref{eq:diff1-1} or \eqref{eq:diff1-2}, including \factanal function, can sometimes suffer from inappropriate solutions. 

Further empirical comparison for each case in Table \ref{tab:numsol} is presented below.   Here, we denote that the $j$th dataset on simulation model S$i$ is S$i$-$j$ ($i=1,2,3; j=1,\dots 10$). 

\subsubsection*{Optimal solution does not exist: S1-2, S1-3, S1-4, S1-5, S2-4}
We first investigate the detailed results in which there are no local maxima. The result (e) in table \ref{tab:numsol} suggests that there exist a number of improper solutions with zero unique variances obtained by \eqref{eq:diff1-2} and \eqref{formula-L}, and some of them provide a positive-definite observed Fisher information matrix (see the result (f)).  However, all of these solutions do not satisfy \eqref{eq:diff0} from the result (g).  Interestingly, the solution computed using \factanal function is close to one of the solutions of \eqref{eq:diff1-2} and \eqref{formula-L} (i.e., the result (f)), suggesting that the \factanal can converge to an inappropriate solution when the optimal solution does not exist. With the \factanal function, the convergence value of the loading matrix constructed by \eqref{eq:diff1-1} does not always satisfy \eqref{eq:diff0} when $\Psi$ is singular; thus, the mixture of update equation \eqref{eq:diff1-1} and Quasi-Newton method with respect to unique variances can lead to the inappropriate solution when the optimal solution does not exist.

The results also give insight into the conventional algorithm of \citet{Jennrich.1969}.  For all simulated datasets, \citeauthor{Jennrich.1969}'s \citeyearpar{Jennrich.1969} algorithm tends to highly depend on initial values, as demonstrated in Figure \ref{fig:discrepancy3}.  Furthermore, one of the unique variances results in an exceptionally negative large value, even if the initial value is close to the solution computed with the {\tt factanal} function.  We observe that the solution obtained by the \citeauthor{Jennrich.1969}'s algorithm \citeyearpar{Jennrich.1969} satisfies \eqref{eq:diff0} but does not satisfy \eqref{eq:diff1-1}.  Recall that \factanal function uses \eqref{eq:diff1-1}; therefore, \factanal function and  \citeauthor{Jennrich.1969}'s \citeyearpar{Jennrich.1969} algorithm result in different solutions.  Furthermore, the observed Fisher information matrix computed with the \citeauthor{Jennrich.1969}'s algorithm \citeyearpar{Jennrich.1969} is not positive definite; one of the eigenvalues of the observed Fisher information matrix is nearly zero, suggesting that the confidence interval is unavailable.  These results are often observed due to the identification problem \citep{Driel.1978,Kano.1998}.  Indeed, the true loading matrix of S1-2, S1-3, S1-4, S1-5 is not unique even if rotational indeterminacy is taken out, because theorem 5.6 in \citet{anderson1956statistical} is not satisfied. % On the other hand, for S2-4, although the true loading matrix is identifiable, there are no optimal solutions.

%\citeauthor{Jennrich.1969}'s \citeyearpar{Jennrich.1969} algorithm tends to highly depend on initial values in our experience, as shown in Figure \ref{fig:discrepancy3}. This is probably because  \citet{Jennrich.1969} use Eq. \eqref{formula-L} that suffers from a number of zero unique variances solutions; indeed, \citet{Jennrich.1969} make use of \eqref{formula-L} to express the loading matrix as a function of unique variances, and find the optimal value of the unique variances using a standard gradient method, such as Newton-Raphson method.  On the other hand, \citeauthor{lawley1971}'s \citeyearpar{lawley1971} algorithm (i.e., the \factanal function) uses Eq. \eqref{eq:diff1-1} that cannot be defined when $\Psi$ is singular, leading to more robust estimates against initial values than \citet{Jennrich.1969}'s algorithm as demonstrated in Figure \ref{fig:discrepancy3}. The \factanal function does not suffer from a large number of zero unique variances solutions in \eqref{formula-L}.  A more detailed discussion will be provided later.

\subsubsection*{Improper solutions: S1-6, S1-7, S2-3, S2-6, S2-7}
We get 5 improper solutions whose unique variances are strictly negative.  These negative solutions are obtained when the proper solutions do not exist.  The number of negative solutions is relatively large for each model from the result (c) of table \ref{tab:numsol}, but most of them do not have a positive-definite observed Fisher information matrix (see result (d)).   As a result, the negative solutions are almost uniquely determined.  An exception is S2-3; there are two solutions that are completely different but have similar likelihood values.

We also find that when the optimal values of unique variances are strictly negative, there often exists a solution to \eqref{eq:diff1-2} and \eqref{formula-L} whose corresponding unique variances are exactly zero and other parameters have similar values to the optimal value. Furthermore, that solution is quite similar to that obtained by the \factanal function.  The corresponding Fisher information matrix is positive-definite.

\subsubsection*{Proper solutions: S1-1, S1-8, S1-9, S1-10, S2-1, S2-2, S2-8, S2-9, S2-10, S3}
We briefly investigate the results where all unique variances are positive.  For all models, we obtain several proper solutions (result (a) in table \ref{tab:numsol}) and only one proper solution has positive-definite observed Fisher information matrix (result (b)).  Furthermore, we obtain no improper solutions whose observed Fisher information matrix is positive (results (d) and (h)), suggesting that the likelihood function has only one local maxima that is proper. 

%\subsubsection*{Identification problems}
\begin{rem}
As shown in Section \ref{mlsgb-2}, our approach computes all candidates for the maximum likelihood solution. If a candidate is a root of a non-zero-dimensional ideal, the candidate may not provide an identifiable model. On the other hand, if the candidate is the root of a zero-dimensional ideal, that candidate does not suffer from identification problems. In our experiments, we often obtain candidates that are roots of a non-zero-dimensional ideal (see \ref{detail-implimentation} for details). Fortunately, however, those candidates do not reach the maximum value of the likelihood.    
\end{rem}

\subsection{Investigation of the best solution obtained by computational algebra and \factanal}
\begin{figure}[!t]
\centering
\includegraphics[width=\textwidth, bb=0 0 972 1040]{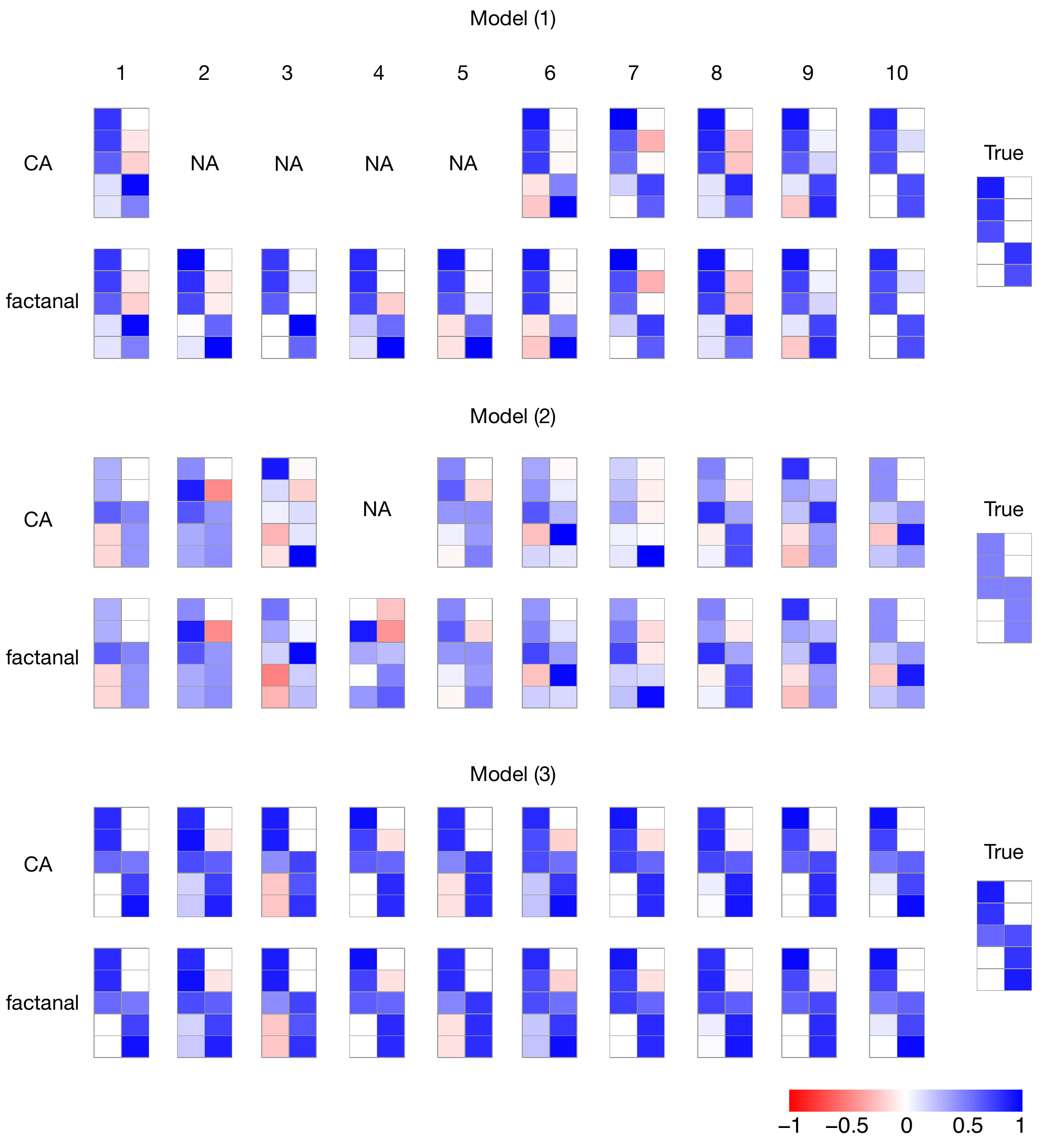}
\caption{Heatmap of the factor loadings obtained by computational algebra and \factanal function.}
\label{fig:heatmap}
\end{figure}
We further compare the characteristics of the loading matrix obtained by computational algebra and {\tt factanal}.  Figure \ref{fig:heatmap} shows the heatmaps of the estimated factor loadings along with the true ones. The loading matrices of computational algebra are unavailable when the optimal solution does not exist (i.e., S1-2, S1-3, S1-4, S1-5, S2-4).  

We first investigate whether the solutions can approximate the true loading matrix.  For S1, the results for the heatmap of the loading matrix in Figure \ref{fig:heatmap} suggest that both computational algebra and \factanal function can well approximate the true loading matrix.  In particular, even if no solution is found with the computational algebra, the \factanal function can approximate the true loadings.  This result is surprising to us because the true loading matrix has an identifiability issue that has been considered as a severe problem \citep{Kano.1998}. In contrast to our expectation, the solutions of S2 are much more unstable than those of S1.  Thus, the small communalities rather than the identification problem would make the estimation difficult in our simulation settings.  For S3, we obtain proper and stable solutions for all settings.

We compare the computational algebra's solutions and {\tt factanal}'s ones. In many cases, these solutions are essentially identical.  In some cases where the computational algebra provides negative unique variances (i.e., S1-7, S2-3, S2-6, S2-7), we find more or less difference between the two methods. This is because the computational algebra provides negative unique variances while the \factanal provides 0.005 due to the restriction of the parameter space.  However, the overall interpretation of these two solutions is essentially identical for S1-7, S2-6 and S2-7.  An exception is S2-3; in this case, the results are entirely different from each other.  Such a difference may be caused by the fact that the \factanal function uses only one initial value. We prepare 100 initial values from $U(0,1)$, perform the \factanal function, and select a solution that maximizes the likelihood function; as a result, the {\tt factanal}'s solution is essentially identical to computational algebra's one.  The result suggests that performing \factanal function with multiple initial values  would be essential when we obtain the improper solutions.

 \subsection{A numerical algorithm to identify the solution pattern}
Because the computational algebra remains a challenge in heavy computation, it would be difficult to apply it to the large model.  Still, our empirical result is helpful for constructing a method to determine whether the solution exists or not with a numerical solution.  Specifically, we obtain the solution by the \citeauthor{Jennrich.1969}'s \citeyearpar{Jennrich.1969} algorithm with multiple initial values.  Subsequently, we employ the same procedure as Algorithm \ref{alg:detectimproper_algsol} based on the numerical solution. The algorithm for identifying the solution pattern with the numerical algorithm is detailed in Algorithm \ref{alg:detectimproper_numsol}.  
 \begin{algorithm}[t]
\caption{A numerical algorithm to identify the cause of improper solutions.}
    \label{alg:detectimproper_numsol}
    \begin{algorithmic}[1] 
    \REQUIRE Dataset \{$\bm{x}_1,\dots,\bm{x}_n$\} or its sample covariance matrix $\bm{S}$, and the number of iterations $T$ to get the maximum likelihood solution.
    \ENSURE Solution pattern.	
    \FOR{$t=1$ to $T$}
	\STATE Obtain a maximum likelihood solution with \citeauthor{Jennrich.1969}'s \citeyearpar{Jennrich.1969} algorithm in \eqref{Ldiff2} using $t$th initial value. 
	\ENDFOR
	\STATE Same as lines 1--10 in Algorithm \ref{alg:detectimproper_algsol}.
    \end{algorithmic}
\end{algorithm}

With algorithm \ref{alg:detectimproper_numsol}, it is possible to further investigate the tendency of solution pattern with a number of simulation iterations.  Table \ref{tab:numsol_sim} shows the distribution of solution pattern with Algorithm \ref{alg:detectimproper_numsol} over 100 simulation runs. The results indicate that only 30\% of S1 provides proper solutions, while S2 has more than double.  Thus, S1 tends to provide more improper solutions than S2. Interestingly, we get several ``improper solution" and ``no solutions" for both S1 and S2, suggesting that various solution patterns are obtained, whether they arise from identifiability issue (S1) or sampling fluctuation (S2).  In other words, it would be difficult to identify the cause of improper solutions by using the solution pattern. However, the distribution of the solution pattern appears to differ among S1 -- S3.  Therefore, it would be essential to find the distribution of solution pattern with a number of simulations (e.g., bootstrap method) for identifying the cause of improper solutions. Such a procedure would be interesting but beyond the scope of this study. Further discussion is described in the next section.
\begin{table}[!t]
  \caption{Solution pattern obtained by Algorithm \ref{alg:detectimproper_numsol} over 100 simulation runs.}
  \label{tab:numsol_sim}
  \begin{center}
  \small
    \begin{tabular}{lcccc}
 \hline
 & Proper solution & Improper solution & No solutions \\ \hline
 S1    & 30    &  36    & 34     \\
 S2    & 66     & 16    & 18    \\
 S3    & 100     & 0    & 0    \\
 \hline
\end{tabular}
\end{center}
\end{table}

 \section{Concluding remarks}
We have proposed an algorithm to get the exact maximum likelihood solutions via computational algebra. The method is based on the algebric equation in \citeauthor{Jennrich.1969}'s \citeyearpar{Jennrich.1969} algorithm.   We obtain all solutions based on Eq. \eqref{Ldiff3}, and select a solution that minimizes the discrepancy function. If the solution's observed Fisher information matrix is positive-definite, that solution is optimal; otherwise, no solution is obtained.  As a result, we can identify the solution pattern of the maximum likelihood estimate with computational algebra.  The conventional analysis may not identify the solution pattern; thus, our proposed approach represents a significant advancement in studying improper solutions in maximum likelihood factor analysis. The improper solutions have been studied for several decades but cannot be resolved from theoretical viewpoints.  We believe that the key to elucidating the improper solution would lie in the advanced pure and applied mathematics, such as computational algebra, singularity theory, and bifurcation theory.  In particular, with advances in computer technology, the computational approaches will play a key role. This study represents a step towards clarifying the improper solutions with advanced applied mathematics. 

For further investigation, we may study transition points among various solution patterns.  For example, we use two sample covariance matrices, S1-1 and S1-2, described in Section \ref{sec:MCS}; S1-1 and S1-2 correspond to positive and no solutions, respectively.  The no solutions of S1-2 is caused by the divergence of $\hat{\psi}_5$. We linearly interpolate these two sample covariance matrices; that is, we construct $S(t) := tS_{1-1} + (1-t)S_{1-2}$ with $t \in [0,1]$.  We note that the maximum likelihood estimation can be performed for any $t \in [0,1]$ because $S(t)$ is positive-definite.  The transition point of positive and negative solutions can be observed by investigating the change in the solution pattern with varying $t$.   

Figure \ref{fig:singular} shows the minimum value of the discrepancy function as a function of $\hat{\psi}_5$ with varying $t$. The results indicate that as $t$ increases, the local minimizer decreases and becomes negative around $t=0.41$.  As $t$ increases further, the discrepancy function becomes flat, and the local minimum disappears at the transition point around $t=0.45$. Beyond this point, there are no global minima.  The discrepancy function varies with a small change in $t$ around 0.4--0.5,  and the change is smooth but drastic.  Such a drastic change would make investigating the improper solutions difficult.  
\begin{figure}[!t]
\centering
\includegraphics[width=\textwidth, bb=0 0 1006 501]{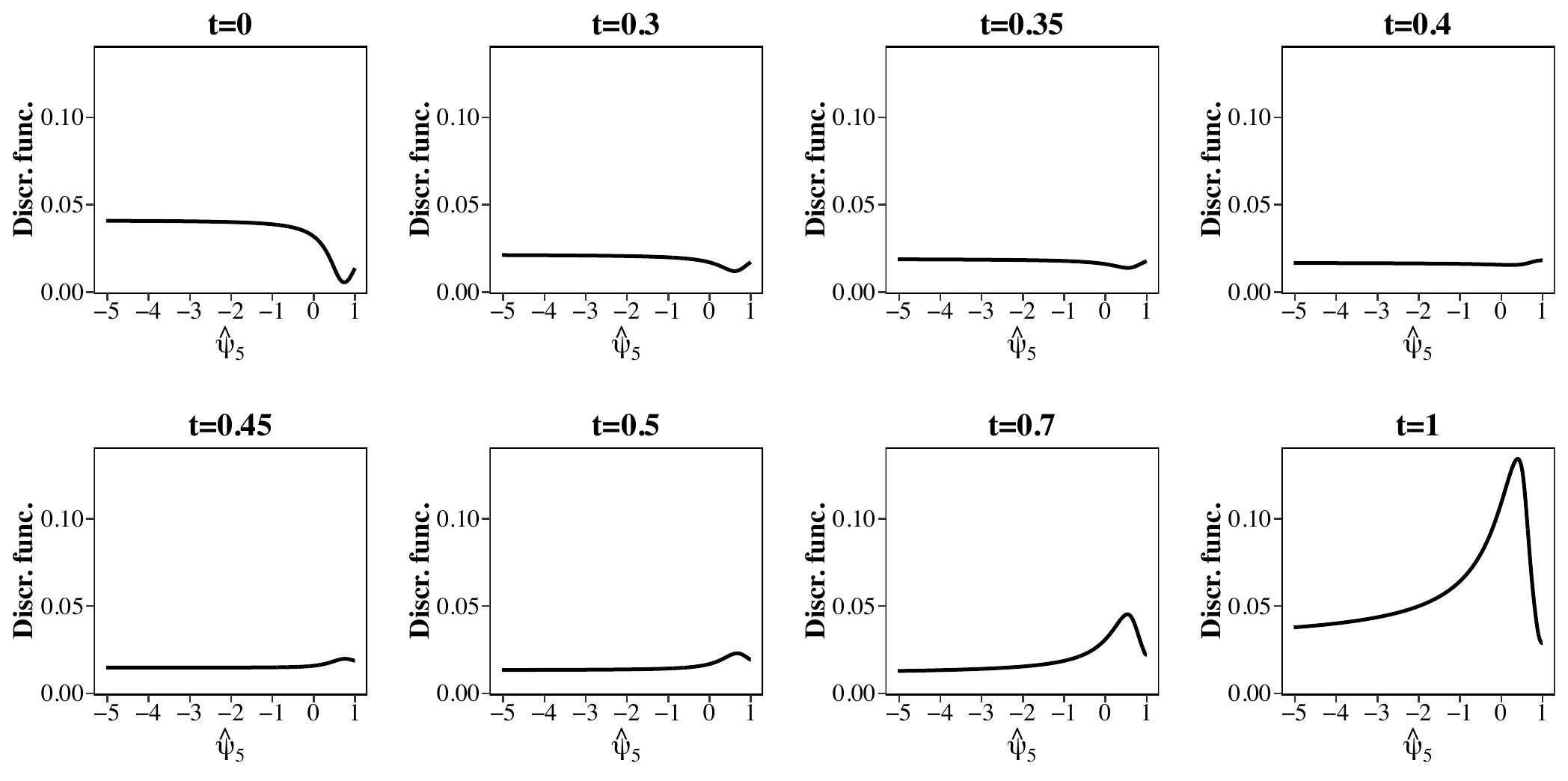}
\caption{Discrepancy function in \eqref{q} as a function of $\hat{\psi}_5$.}
\label{fig:singular}
\end{figure}

For further study, an analysis of geometrical structure based on singularity theory and bifurcation theory would be helpful.  Moreover, in addition to investigating the change in solution pattern between S1-1 and S1-2, it would be interesting to study various change patterns near the transition points with numerous simulations.  These considerations will definitely be significant for further investigating the improper solutions, but they are beyond the scope of this research.  We will take them as future research topics.  

%%%%%%%%%%%%%%%%%%%%%%%%%%%%%%%%%%%%%%%%%%%%%%%%%%%%%%%%%%%%%%%%%%%%%%%%%%%%%%%%%%%%%%%%%%%%%%%%%%%%%%%%%%%%%%%%%%%%%%%%%%%%%%%%%%%%%%%%%%%%%%%%%%%%%%%%%%%%%%%

\bibliographystyle{abbrvnat}
\bibliography{main.bib}

%%%%%%%%%%%%%%%%%%%%%%%%%%%%%%%%%%%%%%%%%%%%%%%%%%%%%%%%%%%%%%%%%%%%%%%%%%%%%%%%%%%%%%%%%%%%%%%%%%%%%%%%%%%%%%%%%%%%%%%%%%%%%%%%%%%%%%%%%%%%%%%%%%%%%%%%%%%%%%%

\appendix
\def\thesection{Appendix \Alph{section}}

\section{Basic Concepts}\label{Fields}

In this section, we will review basic concepts related to computational algebra. In particular, we deal with basic content related to fields, rings and affine varieties along with some concrete examples.  First, we give a definition of fields.

\begin{dfn}\label{dfn:fields}
A field consists of a set $F$ and two binary operations ``$+$'' and ``$\cdot$'' defined on $F$ satisfying the following conditions:
\begin{enumerate}
  \item for any $a,b,c \in F$, $(a+b)+c = a+(b+c)$ and $(a \cdot b) \cdot c = a \cdot (b \cdot c)$ (associativity),
  \item for any $a,b,c \in F$, $a \cdot (b + c) = a \cdot b + a \cdot c$ (distributivity),
  \item for any $a, b \in F$, $a + b = b + a$ and $a \cdot b = b \cdot a$ (commutativity),
  \item for any $a \in F$, there exists $0, 1 \in F$ such that $a + 0 = a \cdot 1 = a$ (identities),
  \item given $a \in F$, there exists $b \in F$ such that $a + b = 0$ (additive inverses),
  \item given $a \in F$, $a \neq 0$, there exists $b \in F$ such that $a \cdot b = 1$ (multiplicative inverses).
\end{enumerate}
\end{dfn}
For example, $\mathbb{Q}$, $\mathbb{R}$ and $\mathbb{C}$ are fields, since they satisfy the following conditions with the sum ``$+$'' and product ``$\cdot$''. On the other hand $\mathbb{Z}$ is not a field, since it does not satisfy the last condition (multiplicative inverses). Indeed, the element $2 \in \mathbb{Z}$ does not have $b \in \mathbb{Z}$ such that $2 \cdot b = 1$.

Next, we define a commutative ring.

\begin{dfn}
A commutative ring consists of a set $R$ and two binary operations ``$+$'' and ``$\cdot$'' defined on $R$ satisfying conditions 1--5 of the Definition \ref{dfn:fields}.
\end{dfn}
As mentioned above, $\mathbb{Z}$ is not a field; however, it is a commutative ring. Moreover $\mathbb{R}[\bm{z}] = \mathbb{R}[z_1, \ldots, z_m]$ is also a commutative ring. In particular, $\mathbb{R}[\bm{z}]$ is known as a polynomial ring.

Then, we introduce the ideals.  In general, the ideals can be defined for any ring. For simplicity, we give a definition of ideals for commutative rings. 
\begin{dfn}
  Let $R$ be a commutative ring. A subset $\mathcal{I} \subset R$ is an ideal if it satifies that
  \begin{enumerate}
  \item $0 \in \mathcal{I}$,
  \item if $a, b \in \mathcal{I}$, then $a + b \in \mathcal{I}$,
  \item if $a \in \mathcal{I}$ and $b \in R$, then $b \cdot a \in \mathcal{I}$.
  \end{enumerate}
\end{dfn}
The definition of ideals is similar to that of linear sub-spaces; both have to be closed under addition and multiplication. However, a difference lies in multiplication; for a linear sub-spaces, we multiply two elements in the field, whereas for ideals, we multiply two elements in the ring. 

In Section \ref{mlsgb-1}, we defined the subset in $\mathbb{R}[\bm{z}]$ as follows:
\begin{align*}
  \langle f_1, \ldots, f_r \rangle = \left\{\sum_{i=1}^r q_i f_i : q_i \in \mathbb{R}[\bm{z}]\right\} \subset \mathbb{R}[\bm{z}],
\end{align*}
where $f_1, \ldots, f_r \in \mathbb{R}[\bm{z}]$. $ \langle f_1, \ldots, f_r \rangle$ is referred to as the ideal generated by $f_1,\ldots, f_r$.  $\langle f_1, \ldots, f_r \rangle$ is similar to the span of a finite number of vectors. In each case, one takes linear combinations, using field coefficients for the span and polynomial coefficients for the ideal. We show that $\langle f_1, \ldots, f_r \rangle$ is an ideal in $\mathbb{R}[\bm{z}]$. 

\begin{lem*}
  $\langle f_1, \ldots, f_r \rangle \subseteq \mathbb{R}[\bm{z}]$ is an ideal.
  \begin{proof}
    Let $\mathcal{J} = \langle f_1, \ldots, f_r \rangle$. If we put $h_1 = \cdots = h_r = 0$, then $0 = \sum_{i=1}^r h_i f_i$. So we have $0 \in \mathcal{J}$.

    Let $h, g \in \mathcal{J}$. The definition of $\mathcal{J}$ implies that $h = \sum_{i=1}^r h_i f_i, g = \sum_{i=1}^r g_i f_i$ for some $h_1,\ldots,h_r, g_1,\ldots,g_r \in \mathbb{R}[\bm{z}]$. Since 
    \begin{align*}
      h + g
      =
      \sum_{i=1}^r h_i f_i + \sum_{i=1}^r g_i f_i
      =
      \sum_{i=1}^r (h_i + g_i) f_i
    \end{align*}
    holds, we have $h + g \in \mathcal{J}$ by $h_1+g_1,\ldots,h_r+g_r \in \mathbb{R}[\bm{z}]$ and the definiton of $\mathcal{J}$.

    Let $h \in \mathcal{J}$ and $c \in \mathbb{R}[\bm{z}]$. Then $h = \sum_{i=1}^r h_i f_i$ holds for some $h_1,\ldots,h_r \in \mathbb{R}[\bm{z}]$. Since
    \begin{align*}
      c h
      =
      c \sum_{i=1}^r h_i f_i
      =
      \sum_{i=1}^r (c h_i) f_i
    \end{align*}
    holds, we have $c h \in \mathcal{J}$ by $c h_1,\ldots,c h_r \in \mathbb{R}[\bm{z}]$ and the definiton of $\mathcal{J}$. Thus we have shown that $\mathcal{J} \subseteq \mathbb{R}[\bm{z}]$ is an ideal.
  \end{proof}
\end{lem*}

Let $U \subseteq W^m$, where $W = \mathbb{R}$ or $\mathbb{C}$. Recall that any affine variety is the solution space in $W^m$ of algebraic equations which has the form \eqref{AEGB-1}. Zariski closures are defined by the following.

\begin{dfn}
  The Zariski closure of $U$, denoted by $\overline{U}$, is the smallest affine variety containing $U$ (in the sense that if $V \subseteq W^m$ is any affine variety containing $U$, then $\overline{U} \subseteq V$).
\end{dfn}

Note that the Zarski closure $\overline{U}$ is the smallest solution space in $W$ containing $U$ among solution spaces of algebraic equations which has the form \eqref{AEGB-1}. We conclude this section with some examples for Zariski closures.

\begin{ex}
  Let $P = \mathbb{V}_{W}(z_1 (z_3^2+1), z_2 (z_3^2+1))$, $Q = \mathbb{V}_{W}(z_3^2+1)$, and $U = P \setminus Q$. Since the subset $U$ is an affine variety, we have $U = \overline{U}$. If $W = \mathbb{R}$, we have $\{ (0,0) \} = \mathbb{V}_{\mathbb{R}}(z_1, z_2) = \overline{U} = U = P$. If $W = \mathbb{C}$, then we have $\{ (0,0) \} = \mathbb{V}_{\mathbb{C}}(z_1, z_2) = \overline{U} = U \subsetneq P$.
\end{ex}

\begin{ex}
  If $P = \mathbb{V}_{W}(z_1 z_2)$, $Q = \mathbb{V}_{W}(z_1, z_2)$, and $U = P \setminus Q$, then $P = \overline{U} \supsetneq  U$.
\end{ex}

\section{Implementation Details}\label{detail-implimentation}

  As a matter of fact, the simulations conducted in Section \ref{sec:MCS} have a very large number of candidates for maximum likelihood solutions. However, since these candidates tend to fall into several types, it is possible to construct more simplified subproblems that further subdivide the subproblems constructed in Algorithm \ref{alg:fukafuka}. In this section, we will show details of these candidates, and present details of the implementation for obtaining the candidates.

Considering a rotational indeterminacy, we suppose that the upper triangular elements of $L$ are fixed by zero in Section \ref{mlsgb-2}. In addition, we consider three models with $p = 5, k=2$ in Section \ref{sec:MCS}. Therefore, we suppose the following form of factor loading matrices in our simulation:
\begin{align*}
  L
  =
  \begin{pmatrix}
    \ell_{11} & 0
    \\
    \ell_{21} & \ell_{22}
    \\
    \vdots & \vdots
    \\
    \ell_{51} & \ell_{52}
  \end{pmatrix}.
\end{align*}

Before discussing details for our implimentation, we show details of solutions to Eqs. \eqref{eq:diff1-2} and \eqref{formula-L}. Table \ref{tab:all-real-number} shows the number of all real solutions to \eqref{eq:diff1-2} and \eqref{formula-L} satisfying that $\Psi$ is not singular, or equivalently
\begin{align*}
  s_{11} - \ell_{11}^2 \neq 0, ~~ s_{22} - (\ell_{21}^2 + \ell_{22}^2) \neq 0, ~~ s_{33} - (\ell_{31}^2 + \ell_{32}^2) \neq 0,
  \\
  s_{44} - (\ell_{41}^2 + \ell_{42}^2) \neq 0, ~~ s_{55} - (\ell_{51}^2 + \ell_{52}^2) \neq 0.
\end{align*}
Note that if $\Psi$ is not singular, the likelihood equation \eqref{eq:diff0} is equivalent to Eqs. \eqref{eq:diff1-2} and \eqref{formula-L}. Therefore, Table \ref{tab:all-real-number} shows the number of all real solutions to \eqref{eq:diff} satisfying that $\Psi$ is not singular. In addition, Table \ref{tab:all-complex-number} shows the number of all complex solutions to \eqref{eq:diff} satisfying that $\Psi$ is not singular. Note that the simulations with infinite complex solutions, namely S1-6, S1-10, S2-10, S3-1, S3-3, S3-4, S3-5, S3-7, form non-zero dimensional ideals. Since we would like to discuss the difficulty of computing all solutions from an algebraic computational aspect, Tables \ref{tab:all-real-number} and \ref{tab:all-complex-number} do not address the identifiability with respect to the sign of column vectors in a loading matrix, unlike Table \ref{tab:numsol}.
Furthermore, Table \ref{tab:numsol} excludes candidates that form a one-factor model, while Tables \ref{tab:all-real-number} and \ref{tab:all-complex-number} do not exclude these candidates.  As shown in the Tables \ref{tab:all-real-number} and \ref{tab:all-complex-number}, Eq. \eqref{eq:diff} has a surprisingly large number of solutions satisfying that $\Psi$ is not singular. In particular, the number of complex solutions is significantly large. Real radical ideals can remove unwanted complex solutions, but they are heavier than Primary ideal decompositions. Therefore, due to such a practical problem, we must skillfully compute real solutions while unfortunately also having unnecessary complex solutions at hand.

\begin{table}[t]  \caption{The number of all real solutions to \eqref{eq:diff} satisfying that $\Psi$ is not singular.}
  \label{tab:all-real-number}
\begin{center}
    \small
    \begin{tabular}{lcccccccccc}
\hline
    Sim\#    & 1     & 2    & 3    & 4     & 5    & 6    & 7    &8    & 9     &10   \\ \hline
    S1 & 29 & 47 & 37 & 33 & 37    & 17 & 53 & 21 & 25 & 13    \\
    S2 & 61 & 25 & 27 & 53 & 33    & 23 & 21 & 27 & 33 & $\infty$ \\
    S3 & 17 & 41 & 29 & 25 & $\infty$ & 33 & 39 & 45 & 49 & 39    \\
\hline
    \end{tabular}
\end{center}
\end{table}

\begin{table}[t]
  \caption{The number of all complex solutions to \eqref{eq:diff} satisfying that $\Psi$ is not singular.}
  \label{tab:all-complex-number}
\begin{center}
    \small
    \begin{tabular}{lcccccccccc}
\hline
    Sim\#    & 1     & 2    & 3    & 4     & 5    & 6    & 7    &8    & 9     &10   \\ \hline
    S1    & 1301 & 1301 & 1301 & 1301 & 1301 & $\infty$ & 1299 & 1301 & 1301 & $\infty$ \\
    S2    & 1301 & 1301 & 1301 & 1301 & 1301 & 1301 & 1301 & 1301 & 1301 & $\infty$ \\
    S3    & $\infty$ & 1301 & $\infty$ & $\infty$ & $\infty$ & 1301 & $\infty$ & 1301 & 1301 & 1301 \\
\hline
    \end{tabular}
\end{center}
\end{table}

Now, let us show the number of all complex solutions to \eqref{eq:diff} satisfying that $\Psi$ is not singular, for each of four spaces satisfying the following conditions; (a) $\ell_{11} = 0$, (b) $\ell_{11} \neq 0$ and $\ell_{i2} = 0$ for some $i = 2, \ldots,5$, (c) $\ell_{11} \neq 0$, $\ell_{i2} \neq 0$ for any $i = 2, \ldots,5$ and $\det L_{J\emptyset} = 0$ for some $J \subset \{ 1, \ldots, 5\}$ satisfying $\#J=3$, (d) $\ell_{11} \neq 0$, $\ell_{i2} \neq 0$ for any $i = 2, \ldots,5$ and $\det L_{J\emptyset} \neq 0$ for any $J \subset \{ 1, \ldots, 5\}$ satisfying $\#J=3$. Table \ref{tab:detnumcmpsol} provides details for Table \ref{tab:all-complex-number} based on these classifications.

\begin{table}[!t]
  \caption{Detailed results for Table \ref{tab:all-complex-number}: the number of all complex solutions to \eqref{eq:diff} satisfying that $\Psi$ is not singular for 4 different parameter spaces.}
  \label{tab:detnumcmpsol}
\begin{center}
  \small
    \begin{tabular}{lrrrrrrrrrl}
\hline
    Sim\#    & 1     & 2    & 3    & 4     & 5    & 6    & 7    &8    & 9     &10   \\ \hline
 \multicolumn{11}{p{12cm}}{(a) All complex solutions to \eqref{eq:diff} satisfying that $\det\Psi\neq0$ and $\ell_{11} = 0$.}   \\ \hline
    S1    & 1 & 1 & 1 & 1 & 1 & $\infty$ & 1 & 1 & 1 & $\infty$ \\
    S2    & 1 & 1 & 1 & 1 & 1 & 1 & 1 & 1 & 1 & $\infty$ \\
    S3    & $\infty$ & 1 & $\infty$ & $\infty$ & $\infty$ & 1 & $\infty$ & 1 & 1 & 1 \\
 \hline\hline
 \multicolumn{11}{p{12cm}}{(b) All complex solutions to \eqref{eq:diff} satisfying that $\det\Psi\neq0$, $\ell_{11} \neq 0$, and $\ell_{i2} = 0$ for some $i = 2, \ldots,5$.}    \\ \hline
    S1    & 212 & 212 & 212 & 212 & 212 & 316 & 210 & 212 & 212 & 346 \\
    S2    & 220 & 212 & 212 & 212 & 212 & 212 & 212 & 212 & 212 & 398 \\
    S3    & 370 & 212 & 288 & 346 & 398 & 212 & 346 & 212 & 212 & 212 \\
 \hline\hline
 \multicolumn{11}{p{12cm}}{(c) All complex solutions to \eqref{eq:diff} satisfying that $\det\Psi\neq0$, $\ell_{11} \neq 0$, $\ell_{i2} \neq 0$ for any $i = 2, \ldots, 5$, and $\det L_{J\emptyset} = 0$ for some $J \subset \{ 1, \ldots, 5\}$ satisfying $\#J=3$.}    \\ \hline
    S1    & 0 & 0 & 0 & 0 & 0 & 48 & 0 & 0 & 0 & 128 \\
    S2    & 0 & 0 & 0 & 0 & 0 & 0 & 0 & 0 & 0 & 128\\
    S3    & 104 & 0 & 112 & 52 & 128 & 0 & 52 & 0 & 0 & 0 \\
 \hline\hline
 \multicolumn{11}{p{12cm}}{(d) All complex solutions to \eqref{eq:diff} satisfying that $\det\Psi\neq0$, $\ell_{11} \neq 0$, $\ell_{i2} \neq 0$ for any $i = 2, \ldots, 5$, and $\det L_{J\emptyset} \neq 0$ for any $J \subset \{ 1, \ldots, 5\}$ satisfying $\#J=3$.}    \\ \hline
    S1    & 1088 & 1088 & 1088 & 1088 & 1088 & 816 & 1088 & 1088 & 1088 & 800 \\
    S2    & 1080 & 1088 & 1088 & 1088 & 1088 & 1088 & 1088 & 1088 & 1088 & 768 \\
    S3    & 784 & 1088 & 896 & 896 & 768 & 1088 & 896 & 1088 & 1088 & 1088 \\
 \hline
\end{tabular}
\end{center}
\end{table}

As can be seen from Table \ref{tab:detnumcmpsol}, the spaces with an infinite number of solutions is limited to the spaces satisfying the condition (a), and we have only a finite number of solutions when the condition (a) is not satisfied. If we use sum ideals and saturation ideals related to the polynomial $\ell_{11}$, we can divide the solution spaces in the spaces which satisfy $\ell_{11} = 0$ or not.

Since the number of solutions satisfying condition (d) is greater than that satisfying conditions (b) and (c), the computation of real solutions satisfying condition (d) seems to require heavy computational loads.  In reality, however, there are several cases where computing real solutions satisfying condition (b) require more heavy computational loads than (c) and (d). The solution space for (c) is more simple than that for (d). For ease of comprehension, we will compare outlines of the algebraic equations satisfying (b) and (d) for a specific dataset, S2-10, which has a large number of solutions satisfying condition (b).  After making the comparison, we present the strategy in our implementation.

Now, let us consider the case of S2-10. A Gr\"obner basis transforms the algebraic equation satisfying $\ell_{22}=0$ into a form consisting of the 15 polynomials $g_{1}^{\mathrm{(b)}},\ldots,g_{15}^{\mathrm{(b)}}$ which have the following leading monomials:
\begin{align*}
  \begin{array}{llll}
  \mathrm{LM}\left(g_{1}^{\mathrm{(b)}}\right) = \ell_{11}, &
  \mathrm{LM}\left(g_{2}^{\mathrm{(b)}}\right) = \ell_{21}, &
  \mathrm{LM}\left(g_{3}^{\mathrm{(b)}}\right) = \ell_{22}, &
  \mathrm{LM}\left(g_{4}^{\mathrm{(b)}}\right) = \ell_{31}^2,
  \\
  \mathrm{LM}\left(g_{5}^{\mathrm{(b)}}\right) = \ell_{31} \ell_{42}, &
  \mathrm{LM}\left(g_{6}^{\mathrm{(b)}}\right) = \ell_{31} \ell_{51}^{74}, &
  \mathrm{LM}\left(g_{7}^{\mathrm{(b)}}\right) = \ell_{31} \ell_{52}, &
  \mathrm{LM}\left(g_{8}^{\mathrm{(b)}}\right) = \ell_{32},
  \\
  \mathrm{LM}\left(g_{9}^{\mathrm{(b)}}\right) = \ell_{41}, &
  \mathrm{LM}\left(g_{10}^{\mathrm{(b)}}\right) = \ell_{42}^2, &
  \mathrm{LM}\left(g_{11}^{\mathrm{(b)}}\right) = \ell_{42} \ell_{51}^{12}, &
  \mathrm{LM}\left(g_{12}^{\mathrm{(b)}}\right) = \ell_{42} \ell_{52},
  \\
  \mathrm{LM}\left(g_{13}^{\mathrm{(b)}}\right) = \ell_{51}^{148}, &
  \mathrm{LM}\left(g_{14}^{\mathrm{(b)}}\right) = \ell_{51}^{2} \ell_{52}, &
  \mathrm{LM}\left(g_{15}^{\mathrm{(b)}}\right) = \ell_{52}^{15},
  \end{array}
\end{align*}
where
\begin{align*}
  \begin{array}{lll}
  g_{1}^{\mathrm{(b)}} \in \mathbb{R}[\ell_{11}, \ell_{21}, \ell_{22}, \ell_{31}, \ell_{32}, \ell_{41}, \ell_{42}, \ell_{51}, \ell_{52}],
  &
  g_{8}^{\mathrm{(b)}} \in \mathbb{R}[\ell_{32}, \ell_{41}, \ell_{42}, \ell_{51}, \ell_{52}],
  &
  g_{15}^{\mathrm{(b)}} \in \mathbb{R}[\ell_{52}].
  \\
  g_{2}^{\mathrm{(b)}} \in \mathbb{R}[\ell_{21}, \ell_{22}, \ell_{31}, \ell_{32}, \ell_{41}, \ell_{42},  \ell_{51}, \ell_{52}],
  &
  g_{9}^{\mathrm{(b)}} \in \mathbb{R}[\ell_{41}, \ell_{42}, \ell_{51}, \ell_{52}],
  \\
  g_{3}^{\mathrm{(b)}} \in \mathbb{R}[\ell_{22}, \ell_{31}, \ell_{32}, \ell_{41}, \ell_{42},  \ell_{51}, \ell_{52}],
  &
  g_{10}^{\mathrm{(b)}} \in \mathbb{R}[\ell_{42}, \ell_{51}, \ell_{52}],
  \\
  g_{4}^{\mathrm{(b)}} \in \mathbb{R}[\ell_{31}, \ell_{32}, \ell_{41}, \ell_{42}, \ell_{51}, \ell_{52}],
  &
  g_{11}^{\mathrm{(b)}} \in \mathbb{R}[\ell_{42}, \ell_{51}, \ell_{52}],
  \\
  g_{5}^{\mathrm{(b)}} \in \mathbb{R}[\ell_{31}, \ell_{32}, \ell_{41}, \ell_{42}, \ell_{51}, \ell_{52}],
  &
  g_{12}^{\mathrm{(b)}} \in \mathbb{R}[\ell_{42}, \ell_{51}, \ell_{52}],
  \\
  g_{6}^{\mathrm{(b)}} \in \mathbb{R}[\ell_{31}, \ell_{32}, \ell_{41}, \ell_{42}, \ell_{51}, \ell_{52}],
  &
  g_{13}^{\mathrm{(b)}} \in \mathbb{R}[\ell_{51}, \ell_{52}],
  \\
  g_{7}^{\mathrm{(b)}} \in \mathbb{R}[\ell_{31}, \ell_{32}, \ell_{41}, \ell_{42}, \ell_{51}, \ell_{52}],
  &
  g_{14}^{\mathrm{(b)}} \in \mathbb{R}[\ell_{51}, \ell_{52}],
  \end{array}
\end{align*}
On the other hand, the algebraic equation of S2-10 satisfying (d), transformed into a form that is easy to solve in the Gr\"obner basis, consists of the 9 polynomials $g_{1}^{\mathrm{(d)}},\ldots,g_{9}^{\mathrm{(d)}}$ which have the following leading monomials:
\begin{align*}
  \begin{array}{lllll}
  \mathrm{LM}\left(g_{1}^{\mathrm{(d)}}\right) = \ell_{11}, &
  \mathrm{LM}\left(g_{2}^{\mathrm{(d)}}\right) = \ell_{21}, &
  \mathrm{LM}\left(g_{3}^{\mathrm{(d)}}\right) = \ell_{22}, &
  \mathrm{LM}\left(g_{4}^{\mathrm{(d)}}\right) = \ell_{31},
  \\
  \mathrm{LM}\left(g_{5}^{\mathrm{(d)}}\right) = \ell_{32}^2, &
  \mathrm{LM}\left(g_{6}^{\mathrm{(d)}}\right) = \ell_{41}, &
  \mathrm{LM}\left(g_{7}^{\mathrm{(d)}}\right) = \ell_{42}, &
  \mathrm{LM}\left(g_{8}^{\mathrm{(d)}}\right) = \ell_{51}^2, &
  \mathrm{LM}\left(g_{9}^{\mathrm{(d)}}\right) = \ell_{52}^{192},
  \end{array}
\end{align*}
where
\begin{align*}
  \begin{array}{lll}
  g_{1}^{\mathrm{(d)}} \in \mathbb{R}[\ell_{11}, \ell_{21}, \ell_{22}, \ell_{31}, \ell_{32}, \ell_{41}, \ell_{42}, \ell_{51}, \ell_{52}],
  &
  g_{5}^{\mathrm{(d)}} \in \mathbb{R}[\ell_{32}, \ell_{41}, \ell_{42}, \ell_{51}, \ell_{52}],
  &
  g_{9}^{\mathrm{(d)}} \in \mathbb{R}[\ell_{52}].
  \\
  g_{2}^{\mathrm{(d)}} \in \mathbb{R}[\ell_{21}, \ell_{22}, \ell_{31}, \ell_{32}, \ell_{41}, \ell_{42},  \ell_{51}, \ell_{52}],
  &
  g_{6}^{\mathrm{(d)}} \in \mathbb{R}[\ell_{41}, \ell_{42}, \ell_{51}, \ell_{52}],
  \\
  g_{3}^{\mathrm{(d)}} \in \mathbb{R}[\ell_{22}, \ell_{31}, \ell_{32}, \ell_{41}, \ell_{42},  \ell_{51}, \ell_{52}],
  &
  g_{7}^{\mathrm{(d)}} \in \mathbb{R}[\ell_{42}, \ell_{51}, \ell_{52}],
  \\
  g_{4}^{\mathrm{(d)}} \in \mathbb{R}[\ell_{31}, \ell_{32}, \ell_{41}, \ell_{42}, \ell_{51}, \ell_{52}],
  &
  g_{8}^{\mathrm{(d)}} \in \mathbb{R}[\ell_{51}, \ell_{52}],
  \end{array}
\end{align*}
It is noteworthy that the leading monomials for (b), namely $\mathrm{LM}(g_i^{\mathrm{(b)}})$, are relatively more complicated than those for (d), namely $\mathrm{LM}(g_i^{\mathrm{(d)}})$, except for the last polynomial $g_{15}^{\mathrm{(b)}}, g_9^{\mathrm{(d)}}$ in each case. All real solutions satisfying condition (d) can be obtained by simple substitutions and simple real number tests once $g_9^{\mathrm{(d)}}$ is factorized, since it has only simple leading monomials except for $g_9^{\mathrm{(d)}}$. On the other hand, it requires heavy computational loads for obtaining real solutions satisfying the condition (b), since there are some leading monomials of high degree, even if $g_{15}^{\mathrm{(b)}}$ is factorized.

Therefore, the following strategy, which decomposes the solution space satisfying condition (b) as delicately as possible, is actually employed.
Let $\mathcal{J} = \langle (\{S-LL^{\top}-\mathrm{diag}(S-LL^{\top})\}S^{-1}L)_{ij} : 1 \leq i \leq 5, 1 \leq j \leq 2 \rangle$. Put
\begin{align*}
  h_1 = s_{11} - \ell_{11}^2, ~~ h_2 = s_{22} - (\ell_{21}^2 + \ell_{22}^2), ~~ \ldots, ~~ h_5 = s_{55} - (\ell_{51}^2 + \ell_{52}^2).
\end{align*}
The above polynoamials are used in Algorithm \ref{alg:fukafuka}. Moreover, we give the polynomials
\begin{align*}
  h_6 &= \ell_{11}, ~~ h_7 = \ell_{22}, ~~ h_8 = \ell_{32}, ~~ h_9 = \ell_{42}, ~~ \ldots, ~~ h_{10} = \ell_{52},
  \\
  h_{11} &= \det L_{\{2,3,4\}\emptyset}, ~~ h_{12} = \det L_{\{2,3,5\}\emptyset}, ~~
  h_{13} = \det L_{\{2,4,5\}\emptyset}, ~~ h_{14} = \det L_{\{ 3,4,5\}\emptyset}.
\end{align*}
First, we construct the sum $\mathcal{I}_0 = \mathcal{J} + \mathcal{K}_{1}$ and the saturation $\mathcal{I}_1 = \mathcal{J}:\mathcal{K}_{1}^{\infty}$. 
Second, we construct the sums $\mathcal{I}_{00} = \mathcal{I}_{0} + \mathcal{K}_{2}, \mathcal{I}_{10} = \mathcal{I}_{1} + \mathcal{K}_{2}$ and saturations $\mathcal{I}_{01} = \mathcal{I}_0:\mathcal{K}_{2}^{\infty}, \mathcal{I}_{11} = \mathcal{I}_1:\mathcal{K}_{2}^{\infty}$. In the same way, for each $b \in \{ 0, 1\}^{i}$, $2 \leq i \leq 13$, we sequentially construct $\mathcal{I}_{b0} = \mathcal{I}_{b} + \mathcal{K}_{i+1}, \mathcal{I}_{b1} = \mathcal{I}_b:\mathcal{K}_{i+1}^{\infty}$. Note that, by Proposition \ref{bunkai}, we have $\mathbb{V}_{W}(\mathcal{J}) = \bigcup_{b \in \{ 0, 1\}^{14}} \mathbb{V}_{W}(\mathcal{I}_{b})$. Finally, we compute all real solutions in $\mathbb{V}_{\mathbb{R}}(\mathcal{I}_b)$ for each $b \in \{0,1\}^{14}$.  
Although there exist cases where such a decomposition provides an empty set of the affine variety for a sum or saturation, the above strategy allowed us to obtain all concrete real solutions.
\end{document}